\documentclass[11pt]{article}

\usepackage{amssymb,amsmath} 
\usepackage{amsthm}
\usepackage{color, graphics}
\usepackage[latin1]{inputenc}
\usepackage[dvips]{graphicx}
\usepackage{natbib}

\title{}
\author{}
\date{}

 \newtheorem{prop}{Proposition}
 \newtheorem{Lemma}{Lemma}
 \newtheorem{thm}{Theorem}
 \newtheorem{cor}{Corollary}
 \newtheorem{df}{Definition}

  \newcommand{\bt}{\begin{thm}}
  \newcommand{\et}{\end{thm}}
  \newcommand{\bdf}{\begin{df}}
  \newcommand{\edf}{\end{df}}
  \newcommand{\bp}{\begin{prop}}
  \newcommand{\ep}{\end{prop}}
  \newcommand{\bc}{\begin{cor}}
  \newcommand{\ec}{\end{cor}}
  \newcommand{\bl}{\begin{Lemma}}
  \newcommand{\el}{\end{Lemma}}

  \newcommand{\be}{\begin{eqnarray*}}
  \newcommand{\ee}{\end{eqnarray*}}
  \newcommand{\ben}{\begin{eqnarray}}
  \newcommand{\een}{\end{eqnarray}}
  \newcommand{\bit}{\begin{itemize}}
  \newcommand{\eit}{\end{itemize}}

\newcommand{\Rex}{\mathbb{R}}

\newcommand{\te}{\theta}
\newcommand{\Te}{\Theta}

\newcommand{\muh}{\hat{\mu}}
\newcommand{\f}{\phi}
\newcommand{\Fi}{\Phi}

\newcommand{\la}{\lambda}

  \newcommand{\gte}{\boldsymbol{\theta}}
  \newcommand{\gTe}{\boldsymbol{\Theta}}
  \newcommand{\gfi}{\boldsymbol{\phi}}
  \newcommand{\gFi}{\boldsymbol{\Phi}}
  \newcommand{\gla}{\boldsymbol{\lambda}}
  \newcommand{\gLa}{\boldsymbol{\Lambda}}
  \newcommand{\gmu}{\boldsymbol{\mu}}

  \newcommand{\gx}{{\bf x}}
  \newcommand{\gX}{{\bf X}}

  \newcommand{\gt}{{\bf t}}
  \newcommand{\gT}{{\bf T}}
  \newcommand{\gG}{{\bf G}}
  \newcommand{\gU}{{\bf U}}
 \newcommand{\gS}{{\bf S}}
  \newcommand{\gs}{{\bf s}}
  
  \newcommand{\gp}{{\bf p}}
 \newcommand{\gZ}{{\bf Z}}
  \newcommand{\gz}{{\bf z}}
    
 \newcommand{\gTqd}{\boldsymbol{T}_{[d]}}
 \newcommand{\gTmqd}{\boldsymbol{T}_{-[d]}}

    \newcommand{\gV}{{\boldsymbol{V}}}
     \newcommand{\gI}{{\boldsymbol{I}}}
    \newcommand{\Hgt}{H_{\boldsymbol{t}}}
    \newcommand{\Hgs}{H_{\boldsymbol{s}}}
    \newcommand{\hgs}{h_{\boldsymbol{s}}}
    \newcommand{\Hgsk}{H_{\boldsymbol{s}_{[k]}}}
    
    \newcommand{\Hgtk}{H_{\boldsymbol{t}_{[k]}, \boldsymbol{t}_{-[d]}}}
    \newcommand{\hgtk}{h_{\boldsymbol{t}_{[k]},\boldsymbol{t}_{-[d]}}}

\newcommand{\s}{\sigma}
\newcommand{\vecx}{x_1,\ldots,x_n}
\newcommand{\vecX}{X_1,\ldots,X_n}

\newcommand{\st}{\vspace{-0.8cm}\begin{flushright} \mbox{$\diamond$} \end{flushright}}

\begin{document}

\begin{center}
{\bf \Large
Fiducial, confidence and objective Bayesian posterior distributions for a multidimensional parameter}

\vspace{4mm}

{\Large Piero Veronese and Eugenio Melilli \\

{\em Bocconi University, Milano, Italy}}


\end{center}

\bigskip

\begin{abstract}

We propose a way to construct fiducial distributions for a multidimensional parameter using a step-by-step conditional procedure related to the inferential importance of the components of the parameter. For discrete models, in which the non-uniqueness of the fiducial distribution is well known, we propose to use the  geometric mean of the \lq\lq extreme cases\rq \rq{} and show its good behavior with respect to the  more traditional arithmetic mean. Connections with the generalized fiducial inference approach developed by Hannig and with confidence distributions are also analyzed. The suggested procedure strongly simplifies when the statistical model belongs to a subclass of the natural exponential family, called conditionally reducible, which includes the multinomial and the negative-multinomial models.
Furthermore, because fiducial inference and objective Bayesian analysis are both attempts to derive distributions for an unknown parameter without any prior information, it is natural to discuss their relationships. In particular, the reference posteriors, which also depend on the importance ordering of the parameters are the natural terms of comparison.  We show that fiducial and reference posterior distributions  coincide in the location-scale models, and we characterize the conditionally reducible natural exponential families for which this happens. The discussion of some  classical examples closes the paper.
 \end{abstract}

\bigskip




\bigskip

\section{Introduction} \label{sec:_intro}

Fiducial distributions, after having been introduced by \citet[1935]{Fisher:1930} and widely discussed (and criticized) in the subsequent years, have been de facto brushed aside for a long time and only recently they have obtained new vitality. The original idea of Fisher was to construct a \emph{distribution} for a parameter which includes all the information given by the data, without resorting to the Bayes theorem. This is obtained by transferring  the randomness from the observed quantity given by the statistical model  to the parameter.
Originally Fisher considered a continuous sufficient statistic $S$ with distribution function $F_{\te}$, depending  on a real parameter $\te$. Let $q_\alpha(\te)$ denote the quantile of order $\alpha$ of $F_{\te}$ and let $s$ be  a realization of $S$. If $q_\alpha(\te)$ is increasing in $\theta$ (i.e., $F_\theta$ is decreasing in $\te$),  the statement  $s<q_\alpha(\te)$ is equivalent to $\theta>q_\alpha^{-1}(s)$ and thus Fisher assumes $q_\alpha^{-1}(s)$ as  the quantile of order $1-\alpha$ of a distribution which he names  \emph{fiducial}.
The set of all quantiles  $q_\alpha^{-1}(s)$, $\alpha \in (0,1)$, establishes the
fiducial distribution function $H_{s}(\te)$ so that
\ben \label{def:original}
H_{s}(\te)=1-F_\te(s) \quad \mbox{and} \quad h_{s}(\te)=\frac{\partial}{\partial \te} H_{s}(\te)= -\frac{\partial}{\partial \te} F_{\te}(s).
\een
Of course $H_s$, and its density $h_s$, must be properly modified if $F_\theta$ is increasing in $\theta$.

\citet[cap.VI]{Fisher:1973}  also provides  some examples of multivariate fiducial distributions obtained by a \lq \lq step-by-step\rq\rq{} procedure, but he never develops a general and rigorous theory. This fact, along with the problem to cover discrete  models, the presence of  some  inconsistencies of the fiducial distribution (e.g. the marginalization paradox, see Dawid \& Stone, 1982), and the difficulties in its interpretation,  gave rise to a quite strong negative attitude towards Fisher proposal.

In the renewed interest for the fiducial approach  a relevant role is played by the \emph{generalized fiducial inference} introduced and developed by \citet[2013]{Hannig:2009}, see also \citet{Hannig:2016} for a review. He provides a formal and mathematically rigorous definition which has  a quite general  applicability. The crucial element of his definition is a data-generating equation $\gX=\gG(\gU,\gte)$, which links the unknown parameter $\gte$ and the observed data $\gX$ through a random element $\gU$ having a known distribution. Roughly speaking,  by shifting  the randomness of $\gU$ from  $\gX$ to $\gte$ (inverting $\gG$ with respect to $\gte$ after having fixed $\gX=\gx$), the distribution given by the statistical model leads to a distribution for the parameter $\gte$. Contrary to the original idea of Fisher, the generalized fiducial distribution is non-unique and Hannig widely discusses this point.
Applications to different statistical models can be found for instance in  \citet{Hannig:2007}, \citet{Hannig:2008} and \citet{Wandler:2012}.

Other recent contributions to the topic of  fiducial distributions are given by \citet{Taraldsen:2013}, \citet{Martin:2013} and Veronese \& Melilli (2015), henceforth  \citet{Veronese:2014}. In this last paper the authors derive  fiducial distributions for  a parameter in a discrete or continuous  real natural exponential family (NEF), and discuss some of their properties with particular emphasis on the frequentist coverage of the fiducial intervals.

In the past fiducial distributions have often been associated with \emph{confidence distributions} even if these latter have a different meaning.
A  modern definition  of confidence distribution is given in \citet{Schweder:2002} and \citet{Singh:2005}, see the book by \citet{Schweder:2016}  for a complete and updated review  on confidence distributions and their connections with fiducial inference.
It is important to emphasize that a confidence distribution  must be regarded as a function of the data with reasonable properties from a purely frequentist point of view.  A confidence distribution is conceptually similar to a point estimator: as there exist several unbiased estimators, several confidence distributions can be provided  for the same parameter and choosing among them can be done resorting to further  optimality criteria. Thus the  confidence distribution theory allows to compare, in a quite general setting, formal distributions for the parameter derived by different statistical procedures.

In this paper we suggest a way to construct a unique distribution for a multidimensional parameter, indexing discrete and continuous models, following a step-by-step procedure similar to that used by \citet{Fisher:1973} in some examples. We call it \emph{fiducial distribution}, but we look at it simply as a distribution on the parameter space in the spirit of the confidence distribution theory.
The key-point of the construction is the procedure by conditioning: the distribution of the data is factorized  as a product of  one-dimensional laws and,  for each of these, the fiducial density for a real parameter component, possibly conditional on other components, is obtained. The joint fiducial density for the parameter is then defined as the product of the (conditional) one-dimensional fiducial densities.
It is well known that  Fisher's fiducial argument presents several drawbacks in higher dimensions,  essentially because one cannot recover the fiducial distribution for a function of the parameters starting from the joint fiducial distribution, see \citet[Ch. 6 and 9]{Schweder:2016}. Our approach, when it can be applied, presents the advantage to construct sequentially the fiducial distribution directly on the parameters of interest and different fiducial distributions can be obtained focusing on  different parameters of interest.  Also, it should be noticed that a general definition  of confidence distribution for a  multidimensional parameter does not exist and more attention is given to the construction of approximate confidence curves for specific nested families of regions, see  \citet[Ch. 9 and Sec. 15.4]{Schweder:2016}.

Interestingly, our joint fiducial distribution coincides in many cases with the Bayesian posterior obtained using the reference prior.
This fact motivates the second goal of the paper: to investigate  the relationships between the objective Bayesian posteriors and the suggested fiducial distributions.
Objective Bayesian analysis, see e.g. \citet{Berger:2006}, essentially studies how to perform a \emph{good}  Bayesian inference, especially for moderate sample size, when one is unwilling or unable to assess a subjective prior. Under this approach, the prior distribution is derived directly from the model and thus  it is labeled as \emph{objective}. The reference prior, introduced by \citet{Bernardo:1979} and developed by \citet{Berger:1992}, is the most successful default prior proposed in the literature.
For a multidimensional parameter the reference prior depends on the grouping and ordering of its  components and, in general, no longer coincides with the Jeffreys prior. This is  the reference prior only for a real parameter and  it is unsatisfactory otherwise, as well known.

\citet{Lindley:1958} was the first to discuss the connections between fiducial and posterior distributions for a real parameter, when a real continuous  sufficient statistic exists.
\citet{Veronese:2014} extend this result to real discrete NEFs,
characterizing  all families admitting a  \emph{fiducial prior}, i.e. a prior leading to a posterior coinciding with the  fiducial distribution.
This prior is strictly related to  the Jeffreys prior.
We show here that when the parameter is multidimensional this relationship no longer holds and a new one is established with the reference prior. In particular we prove results for location-scale parameter models and \emph{conditionally reducible} NEFs, a subclass of NEFs defined in \citet{Consonni:2001}.

The paper is structured  as follows.
Section \ref{sec:preli}  reviews some basic facts  on fiducial and confidence distributions for real  NEFs and on generalized fiducial distributions.
 The proposal for constructing a step-by-step multivariate fiducial distribution  is presented in Section \ref{sec:step-by-step fiducial}, which also discusses: the relationships with confidence distributions   (Section \ref{sec:relat_with_conf}), the use of the geometric mean
of fiducial densities for solving the non-uniqueness problem in discrete models (Section \ref{sec:discrete case}),  the connections with the generalized fiducial inference and the consistency with the sufficiency principle (Section \ref{sec:generalized fiducial}).
Section \ref{sec:cr-NEF} studies the fiducial distributions for conditionally reducible NEFs and provides their explicit expression for a particular subclass which includes the  multinomial and the negative-multinomial model.
Section \ref{sec:connection reference} analyzes  the relationships between the fiducial distributions and the reference posteriors, in particular for location-scale parameter models (Section \ref{sec:location-scale}) and NEFs (Section \ref{sec:bayesian-NEFs}), characterizing those which  admit  the fiducial prior.
  Section \ref{sec:further examples} discusses  further examples in which fiducial and reference posteriors coincide. Section  \ref{sec:conclusions} concludes the paper presenting some possible asymptotic extensions.
Finally, Appendix A1 collects some  useful technical results on conditionally reducible NEFs, while Appendix A2 includes the proofs of all the results stated in the paper.

\section{Preliminary results}\label{sec:preli}

 The modern definition of  \emph{confidence distribution} for a real parameter $\phi$ of interest, see \citet[2016]{Schweder:2002} and \citet{Singh:2005}, can be formulated as follows:
\bdf \label{CD}
Let $\{F_{\phi,\gla},\phi \in \Phi \subseteq \Rex, \gla \in \gLa\}$ be a parametric model for data $\gX \in \cal X$; here $\phi$ is the parameter of interest and $\gla$ is a nuisance parameter. A function $C:{\cal X} \times \Phi \rightarrow \Rex$ is a \emph{confidence distribution} for $\phi$ if $C(\gx,\cdot)$ is a distribution function for each $\gx \in \cal X$ and $C(\gX,\phi_0)$ has a uniform distribution in $(0,1)$ under $F_{\phi_0,\gla_0}$,
where $(\phi_0,\gla_0)$ is the true parameter value.
\edf
The relevant requirement in the previous definition  is the uniformity of the distribution, which ensures the correct coverage of the confidence intervals. As seen in Section \ref{sec:_intro}, the confidence distribution theory must be placed in  a purely frequentist context and allows to compare distributions on the parameter space, obtained using different approaches. Finally, the definition of confidence distribution can be generalized by requiring that the uniformity assumption holds only asymptotically.

Strictly linked to the notion of confidence distribution is that of \emph{confidence curve}, defined, for each observed $\gX=\gx$,  as the function  $\phi \rightarrow cc(\phi)=|1-2C(\gx,\phi)|$; see \citet{Schweder:2016}. This function gives the extremes of equal-tail confidence intervals  for any level $1-\alpha$,  allowing a fast and clear comparison of confidence distributions with respect to their interval length.
When the parameter of interest is multidimensional,  how to extend the definitions of  confidence distribution and confidence curve is much less clear and various proposals have been made, see \citet[2016]{Schweder:2002} and  \citet{Singh:2005}.

As detailed in Section \ref{sec:_intro}, \citet{Hannig:2009}  has proposed the notion of \emph{generalized fiducial distribution}, which is based on a data-generating equation  $\gX=\gG(\gte,\gU)$. Because several functions $\gG$ can  generate the same statistical model, and not all the resulting fiducial distributions  are reasonable in terms of properties or computational tractability, \citet[Sec. 5]{Hannig:2013} gives some hints on the choice of a \emph{default} function $\gG$. In particular, if $\gX=(X_1, \dots,X_n)$  is an independent identically distributed (i.i.d.) random sample from an (absolutely) continuous distribution function $F_{\gte}$, with density $f_{\gte}$, $\gte \in \Rex^d$, he suggests to use $X_i=F^{-1}_{\gte}(U_i)$, $i=1,\dots,n$, where $U_i$ are i.i.d. uniform random variables on $(0,1)$ and $F^{-1}_{\gte}$ is the inverse (or generalized inverse) of $F_{\gte}$.  If other regularity assumptions are satisfied, the generalized fiducial distribution  for $\gte$ can be written as
\ben \label{r-hannig}
r(\gte)=\frac{f_{\gte}(\gx) J(\gx,\gte)}{\int_{\gTe}f_{\gte}(\gx) J(\gx,\gte) d \gte},
\een
where the expression of  $J(\gx,\gte)$,  given in \citet[formula (3.7)]{Hannig:2013}, is
\ben \label{Hannig-J}
J(\gx,\gte)= \sum_{\{(i_1, \ldots, i_d):1\leq i_1< \cdots < i_d \leq n\}} \left|\frac{\det\left(\frac{d}{d\gte}(F_{\gte}(x_{i_1}),\dots,F_{\gte}(x_{i_d}))\right)}{\prod_{j=1}^d f_{\gte}(x_{i_j})}\right|.
\een
In \eqref{Hannig-J}  the numerator of the ratio is the determinant of the matrix whose $kj$-entry is $\partial F_{\gte}(x_{i_j})/\partial \te_k$.
This procedure leads to the Fisher definition  of fiducial density \eqref{def:original} when $n=d=1$.

\citet[Example 4]{Hannig:2013} explicitly recognizes the advise of \citet{Wilkinson:1977} that the choice of a fiducial distribution should depend on the parameter of interest and uses the well known example of $d$ independent normal distributions N$(\mu_i,1)$, in which the parameter of interest is $\te=(\sum_{i=1}^d \mu_i^2)^{1/2}$. He shows that the default data-generating equations $X_i=\mu_i+U_i$, $i=1,\dots,d$, lead to a fiducial distribution  which has good frequentist properties for inference on the $\mu$'s, but very bad ones when the interest is on $\te$, as already recognized by \citet{Stein:1959}. Thus Hannig suggests an ad hoc alternative  equation, which leads to a better solution.
Notice that our general procedure, suggested in the next section, constructs a fiducial distribution  starting directly from the parameter of interest and do not required the choice a priori of a data generating function.

Fiducial distributions and their properties, with particular emphasis on the frequentist coverage of the fiducial intervals, for a  discrete or a continuous real regular NEF, are discussed in \citet{Veronese:2014}.
More specifically,  consider the sufficient statistic $S$ associated with a sample of size $n$ and denote by ${\cal S}$ its support. Let    $F_\theta (s)$ be the distribution function of $S$ and  $p_{\theta}(s) =  \exp\left\{ \theta s  - nM(\theta) \right \} $ the corresponding density (with respect to a measure $\nu$).
 Let $a=\inf {\cal S}$, $b=\sup {\cal S}$ and
define ${\cal S}^*=[a,b)$ if $\nu(a)>0$, otherwise
${\cal S}^*=(a,b)$.   Then,  for $s \in {\cal S}^*$, \citet{Petrone:2010} have proved that
 \ben \label{Hns}
H_s(\te)= \left\{ \begin{array}{ll} 0 & \hspace{5mm} \te \leq \inf \Theta\\
1-F_{\te}(s) & \hspace{5mm}   \inf \Theta < \te <  \sup \Theta\\ 1 &  \hspace{5mm}   \te \geq \sup \Theta
\end{array} \right .
\een
  is  a  fiducial distribution function  for the natural parameter $\te$. It follows that the fiducial  density of $\te$ is
 \ben \label{hns}
  h_s(\te) =  \frac{\partial}{\partial\te} H_s(\te) = -\frac{\partial}{\partial \te} F_{\te}(s)  = \int_{(-\infty,s]}(nM'(\te)-t)p_{\te}(t)
   d\nu(t).
 \een
It is important to underline, and  simple to verify, that the distribution function $H_s$ is also a  confidence distribution (only asymptotically, in the discrete case), according to Definition \ref{CD}.

Notice that, for discrete NEFs,  $F_{\theta}(s)=\mbox{Pr}_\te\{S \leq s\}$ and  Pr$_\te\{S<s\}$  do not coincide and thus, besides   $H_s$
in (\ref{Hns}), one could define a
\emph{left fiducial distribution} as
\ben \label{left-fid}
H_s^\ell(\te)=1-\mbox{Pr}_\te\{S<s\}.
\een
 For convenience, sometimes $H_s$ will be called \emph{right fiducial distribution}.  A standard way to overcome this non-uniqueness is referring to the half-correction device \citep[see][pag. 62]{Schweder:2016} which amounts to consider the mixture  $H_s^A(\te)=(H_s(\te) +H_s^\ell(\te))/2=\mbox{Pr}_{\te}\{S>s\}+\mbox{ Pr}_{\te}\{S=s\}/2$, whose density is  the arithmetic mean of $h_s(\te)$ and $h_s^\ell(\te)$.
Instead, we will suggest  to average  $h_s$ and $h_s^\ell$ using their \emph{geometric mean} $h_s^G$ (suitably normalized) and show that it presents better properties than  $h_s^A$  (Section \ref{sec:discrete case}) and   a more direct connection with objective Bayesian inference (Section \ref{sec:bayesian-NEFs}), even if, operationally, the difference is usually  not particularly big.

Table \ref{tab_nef-qvf} provides the fiducial distributions obtained in \citet{Veronese:2014}  for some important discrete and continuous NEFs, which will be used in the forthcoming examples. It also establishes the abbreviations used in the paper for the standard distributions.

\begin{table}[!hbp]
 \caption{ \emph{Fiducial distributions for some real NEFs\label{tab_nef-qvf}}}
\begin{center}
\footnotesize{
{%
\begin{tabular}{lll}
  \hline
  & Sufficient        &   Fiducial                         \\
  &   statistic               &          distributions                  \\ \hline
  N$(\mu, \sigma^2)$  & $S=\sum_i X_i$ & $H_s(\mu)$ : N$(s/n,\sigma^2/n)$  \\
  ($\sigma^2$ known) &  &   \\
  \\
  N$(\mu, \sigma^2)$  & $S=\sum_i (X_i-\mu)^2$ &  $H_s(\sigma^2)$: In-Ga$(n/2,s/2)$ \\
  ($\mu$ known) & &  \\
  \\
  Ga($\alpha, \lambda$) & $S=\sum_i X_i$ & $H_s(\lambda)$ : Ga($n \alpha,s$) \\
  ($\alpha$ known) & &      \\ \\
  Pa$(\lambda, x_0)$  & $S=\sum_i \log (X_i/x_0)$ &  $H_s(\lambda)$ : Ga$(n,s)$ \\
  ($x_0$ known) & &   \\
  \\
  We$(\lambda, c)$  & $S=\sum_i X_i^c$ &  $H_s(\lambda)$ : Ga$(n,s)$   \\
  ($c$ known) & &   \\
  \\
  Bi($m,p$) &  $S=\sum_i X_i$ &  $H_s(p)$ : Be$(s+1,nm-s)$   \\
  ($m$ known) &  &   $H_s^\ell(p)$ : Be$(s,nm-s+1)$    \\
            &    &  $H_s^G(p)$ : Be$(s+1/2,nm-s+1/2)$     \\
  \\
  Po($\mu$) & $S=\sum_i X_i$ &  $H_s(\mu)$ : Ga($s+1,n$)   \\
   &    &  $H_s^\ell(\mu)$ : Ga$(s,n)$           \\
    &    &  $H_s^G(\mu)$ : Ga$(s+1/2,n)$         \\
  \\
  Ne-Bi($m,p$)  &  $S=\sum_i X_i$   & $H_s(p)$ : Be($nm,s+1$)  \\
  ($m$ known) & & $H_s^\ell (p)$ : Be($nm,s$)        \\
             &  &      $H_s^G(p)$ : Be($nm,s+1/2$)          \\ \hline
  \end{tabular}}}
  \end{center}
\footnotesize{\emph{The following notations are used:
 $\mbox{Ga}(\alpha,\lambda)$ for a  gamma distribution with shape $\alpha$ and mean  $\alpha/\lambda$;
 $\mbox{In-Ga}(\alpha,\lambda)$ for an inverse-gamma distribution (if $X \sim \mbox{Ga}(\alpha,\lambda)$ then $1/X \sim \mbox{In-Ga}(\alpha,\lambda))$;
 $\mbox{Be}(\alpha,\beta)$ for a beta distribution with parameters $\alpha$ and $\beta$;
 $\mbox{Bi}(m,p)$ for a  binomial distribution with $m$ trials and success probability $p$;
 $\mbox{Ne-Bi}(m,p)$ for a  negative-binomial with  $m$ successes and success probability $p$;
 $\mbox{Po}(\mu)$ for the Poisson distribuition with mean $\mu$;
 $\mbox{Pa}(\lambda, x_0)$ for a Pareto distribution with density  $\la x_0^\lambda x^{-\la-1}$, $x>x_0>0$, $\la>0$;
 $\mbox{We}(\lambda, c)$ for a  Weibull distribution with density $c\la x^{c-1}\exp(-\la x^c)$,  $x, \la,c >0$}.}
\end{table}

\section{Fiducial distributions for multidimensional parameters} \label{sec:step-by-step fiducial}

 A  \emph{natural} way to construct a suitable fiducial distribution for a multidimensional  parameter is to follow the step-by-step procedure used by \citet{Fisher:1973} in some examples. The key-point of our proposal stems on the factorization of the sampling distribution as a product of one-dimensional conditional laws. For each of these the fiducial density for a real component of the parameter, possibly conditional on other components, is defined.
It is well known that different factorizations of sampling distributions can produce different joint fiducial distributions, see e.g. \citet{Dempster:1963}. However,  we do not consider this aspect a drawback of the procedure if it is linked to the inferential importance ordering  of the parameter components implied by the  factorization. For example, if a parameter $\gte \in \Rex^2$ is transformed in such a way that $\f$ is the parameter of interest and $\la$ the nuisance, the obvious ordering is  $(\f,\la)$ and a suitable  factorization must be defined accordingly,  see Example 4 (ctd.) in this section for an illustration. The crucial role played by the ordering of the parameters accordingly to their inferential importance is widely acknowledged by objective Bayesian inference, in which reference priors are different for different orderings, see Section \ref{sec:connection reference}.

In order to construct a fiducial distribution, we consider two basic transformations: one involving  the sample data $\gX=(X_1,\dots,X_n)$, having a distribution parameterized by $\gte=(\te_1,\dots,\te_d)$,  $d\leq n$, and one involving $\gte$.
Given $\gX$, consider a statistic $\gT=(T_1, \dots, T_m)$, $d\leq m\leq n$, with density  $p_{\gte}(\gt)$, which summarizes $\gX$ without losing information on $\gte$. $\gT$ can be a sufficient statistic or a one-to-one transformation of $\gX$. Split $\gT$ in  $(\gT_{[d]}, \gT_{-[d]})$, where $\gT_{[d]}=(T_1, \dots, T_d)$ and  $\gT_{-[d]}= (T_{d+1}, \dots, T_m)$,  and suppose that $\gT_{-[d]}$ is ancillary for $\gte$. As a consequence $p_{\gte}(\gt)= p_{\gte}(\gt_{[d]}|\gt_{-[d]})p(\gt_{-[d]})$ and all the information on $\gte$ provided by $\gX$ are included in the conditional distribution of $\gT_{[d]}$ given $\gT_{-[d]}$.

Assume now that there exists a one-to-one smooth reparameterization from $\gte$ to $\gfi$, with $\f_1, \ldots, \f_d$ ordered with respect to  their importance, such that
\ben \label{general-conditional}
p_{\gfi}(\gt_{[d]}|\gt_{-[d]})=\prod_{k=1}^d p_{\f_{d-k+1}}(t_k|\gt_{[k-1]}, \gt_{-[d]};\gfi_{[d-k]}).
\een
The density $p_{\f_{d-k+1}}(t_k|\gt_{[k-1]}, \gt_{-[d]};\gfi_{[d-k]})$, with  the corresponding distribution function $F_{\f_{d-k+1}}(t_k|\gt_{[k-1]},\gt_{-[d]};\gfi_{[d-k]})$, must be interpreted as the conditional distribution of $T_k$ given $(\gT_{[k-1]}=\gt_{[k-1]}, \gT_{-[d]}=\gt_{-[d]})$, parameterized by $\f_{d-k+1}$, assuming $\gfi_{[d-k]}$  known.
In the following,  we will always assume that all the one-dimensional conditional distribution functions  $F_{\f_j}$'s involved in the analysis  are monotone and differentiable  in $\f_j$ and have limits $0$ and $1$ when $\f_j$ tends to the boundaries of its domain.
Notice that this is always true if $F_{\f_j}$ belongs to a NEF, see \eqref{Hns}.
Under these assumptions, the joint fiducial density of $\gfi$ is obtained as
\ben \label{fid-general}
h_{\gt}(\gfi)=\prod_{k=1}^d h_{\gt_{[k]}, \gt_{-[d]}}(\f_{d-k+1}|\gfi_{[d-k]}),
\een
 and
\ben \label{fid-general-2}
h_{\gt_{[k]}, \gt_{-[d]}}(\f_{d-k+1}|\gfi_{[d-k]})= \left|\frac{\partial}{\partial \f_{d-k+1}} F_{\f_{d-k+1}}(t_k|\gt_{[k-1]},\gt_{-[d]}; \gfi_{[d-k]})\right|.
\een

Several applications of this procedure to well known models will be provided in Section \ref{sec:further examples}.
Here we illustrate some interesting features of the fiducial distribution \eqref{fid-general}.

\medskip

\noindent
i)
The existence of an ancillary statistic is not necessary if there exists a sufficient statistic with the same dimension of the parameter $(m=d)$. An important case is $m=d=1$ so that formula  \eqref{fid-general} and \eqref{fid-general-2} reduce to
$h_t(\f)=\left|\partial F_{\f}(t)/\partial \f\right|$,
the original formula suggested by \citet{Fisher:1930}. \\

\medskip

\noindent
ii)
 If one is only interested in $\f_1$, it follows from \eqref{general-conditional} that it is enough to consider
 \be
 h_{\gt}(\f_{1})= \left|\frac{\partial}{\partial \f_1} F_{\f_{1}}(t_d|\gt_{[d-1]},\gt_{-[d]})\right|,
 \ee
 which, depending on all observations, does not lose any sample information.
 A typical choice for  $T_{d}$ is given by the maximum likelihood estimator $\widehat{\f}_1$ of $\f_1$ and thus, when $\widehat{\f}_1$ is not sufficient, we have to consider the distribution of $\widehat{\f}_1$  given the ancillary statistic $\gt_{-[d]}$.
 Similarly, if one is interested in $\f_1, \f_2$, it is enough to consider $ h_{\gt}(\f_{1})\cdot h_{\gt_{[d-1]}, \gt_{-[d]}}(\f_{2}|\f_{1})$, and so on.

\medskip

\noindent
iii)
 When an ancillary statistic $\gT_{-[d]}$ is needed, the fiducial distribution \eqref{fid-general} is invariant with respect to any one-to-one transformation of $\gT_{-[d]}$. All the sampling distributions are conditional on it and thus any transformation establishes the same constraints; see Section \ref{sec:location-scale} for an example.

\medskip

\noindent
iv)
The construction by successive conditioning makes the fiducial distribution invariant under the so called one-to-one \emph{lower triangular transformation} of $\gT_{[d]}$, for fixed $\gT_{-[d]}$. More precisely,
we consider a transformation $\gT^*=(\gT^*_{[d]}, \gT_{-[d]})$
   such that $T^*_k=g_k(\gT_{[k]},\gT_{-[d]})$, for  $k=1,\dots,d$. To see this, assuming for instance
  $t^*_{k}=g_k(\gt_{[k]}, \gt_{-[d]})$ increasing in $t_k$, it is sufficient to show that
\be
\lefteqn{\mbox{Pr}_{\f_{d-k+1}}(T^*_k\leq t^*_k \mid\gT^*_{[k-1]}=\gt^*_{[k-1]},\gT_{-[d]}=\gt_{-[d]}; \gfi_{[d-k]})}\\
& & = \mbox{Pr}_{\f_{d-k+1}}(g_k(\gT_{[k]},\gT_{-[d]}) \leq g_k(\gt_{[k]},\gt_{-[d]})\mid\gT^*_{[k-1]}=\gt^*_{[k-1]},\gT_{-[d]}=\gt_{-[d]}; \gfi_{[d-k]})\\
& & = \mbox{Pr}_{\f_{d-k+1}}(T_{k} \leq t_{k}\mid\gT_{[k-1]}=\gt_{[k-1]},\gT_{-[d]}=\gt_{-[d]}; \gfi_{[d-k]}).
\ee
It follows immediately that $\gT$ and $\gT^*$ lead to the same fiducial distribution.

\medskip

\noindent
v) If $(\gT_{[k-1]}, \gT_{-[d]})$ is sufficient for $\gfi_{[d-k]}$, for each $k$, then the conditional distribution of   $T_k$ given $(\gT_{[k-1]}=\gt_{[k-1]}, \gT_{-[d]}=\gt_{-[d]})$  does not depend on
  $\gfi_{[d-k]}$ and the fiducial distribution \eqref{fid-general} becomes the product of the \lq\lq marginal\rq\rq{} fiducial distributions of the $\f_{k}$'s. As a consequence,  \eqref{fid-general-2} can be used alone to make inference on $\f_{d-k+1}$ and the fiducial distribution does not depend on the inferential ordering of the parameters. An important case in which this happens  will be discussed in Section \ref{sec:cr-NEF}.

\medskip

We close this section establishing the invariance property of the fiducial distribution $h_{\gt}(\gfi)$ under a lower triangular transformation, i.e. a transformation  from $\gfi$ to $\gla=(\la_1, \ldots, \la_d)$, say, which maintains the same decreasing ordering  of importance in the components of the two vectors.

\begin{prop} \label{prop:invariance}
If $\gfi=\gfi(\gla)$ is a one-to-one lower triangular continuously differentiable function from $\gLa$ to $\gFi$, then the fiducial distribution $h_\gt^{\gfi}(\gfi)$, obtained applying \eqref{fid-general} to the  model $p_{\gfi}(\gt)$, and the fiducial distribution $h_\gt^{\gla}(\gla)$, obtained applying \eqref{fid-general} to the  model  $p_{\gla}(\gt)=p_{\gfi(\gla)}(\gt)$, are such that, for each measurable  $A \subset \gFi$,
\ben \label{inv-fid}
\int_{A}h_\gt^{\gfi}(\gfi)d \gfi= \int_{\gla^{-1}(A)}h_\gt^{\gla}(\gla) d\gla.
\een
\end{prop}

\subsection{Relationships with confidence distributions}\label{sec:relat_with_conf}

Given a real NEF, $H_{s}(\te)$ in \eqref{Hns} is an exact or approximate confidence distribution if the observations are continuous or discrete, respectively. It is possible to verify that the same is true for the marginal fiducial distribution of the main parameter of interest  $\phi_1$ in the more general definition (\ref{fid-general}). Indeed, the  distribution function of $\phi_1$ is
$H_{\gt}(\f_{1})= 1- F_{\f_{1}}(t_d|\gt_{[d-1]},\gt_{-[d]})$,
so that the first requirement  in Definition \ref{CD} is clearly satisfied, thanks to the assumption on the distribution function given after formula  (\ref{general-conditional}). For what concerns the  uniformity condition, assuming that $F_{\phi_1}$  is decreasing in $\phi_1$ (if it is increasing,  replace $1-F_{\phi_1}$ with $F_{\phi_1}$), we have, for $u \in (0,1)$ and arbitrary $\gfi$,
\be
\mbox{Pr}_{\gfi} \left(H_{\gt_{[d]}, \gt_{-[d]}}(\f_{1}) \leq u \right) = 1-\mbox{Pr}_{\gfi} \left(F_{\f_{1}}(t_d|\gT_{[d-1]},\gT_{-[d]}) < 1-u\right)=
\ee
\be
=1-\int \mbox{Pr}_{\gfi}\left\{ F_{\f_{1}}(t_d|\gt_{[d-1]},\gt_{-[d]}) < 1-u\right\} dF_{\gfi}(\gt_{[d-1]},\gt_{-[d]})=u
\ee
because, by construction,  the integrand  is equal to $1-u$ for all fixed $(\gt_{[d-1]},\gt_{-[d]})$.

\subsection{The discrete case: the geometric mean of the left and right fiducial densities} \label{sec:discrete case}

As mentioned in Section \ref{sec:preli}, for a discrete statistic $S$ with distribution depending on a real
parameter $\te$, we suggest to use the  geometric mean of the right and left fiducial densities,  $h_{s}^G(\theta)=c^{-1}(h_{s}(\theta) h_{s}^\ell (\theta))^{1/2}$, where $c$ is the normalizing constant,  instead of their arithmetic mean $h_{s}^A(\theta)$.

A first justification of the use of the geometric mean of densities  is suggested by \citet{Berger:2015} who  mention its property to be  the density \lq \lq closest\rq\rq{}  to $h_s$ and $h_s^\ell$ with respect to the the Kullback-Leibler divergence, as specified in the following proposition. We give a simple proof of this fact, without resorting to the calculus of variations. Recall that,
given two densities $p$ and $q$,  having the same support and the same dominating measure $\nu$, the Kullback-Leibler divergence of $p$ from $q$ is defined as $KL(q|p)= \int q(x) \log (q(x)/p(x))d\nu(x)$.

\bp \label{prop: KL-geom-mean}
Consider two  densities $p_1$ and $p_2$ with the same support. The density $q$ which
minimizes $KL(q|p_1)+KL(q|p_2)$ is  given by
$q=p^G\propto (p_1p_2)^{1/2}$, which is the (normalized)  geometric mean of $p_1$ and $p_2$.
\ep

Furthermore,  \citet{Krishnamoorthy:2010} observe that a distribution for $\theta$, whose aim is to give a synthesis of two fiducial distributions,   should ``stochastically'' lie between them. In our setting, the extreme distributions are  $H_{s}$ and $H_{s}^\ell$.  This property is surely satisfied by  the arithmetic mean, because  $H_{s} (\theta) < H^A_{s} (\theta) < H_{s}^\ell(\theta)$ uniformly with respect to $\theta$, for each $s$ belonging to the set ${\cal S}_0$ for which both $H_{s}$ and $H_{s}^\ell$ can be defined. The same inequalities are  true for $H^G_s$ under mild assumptions. As usual, here we assume that  $H_s(\theta)$ is defined as $1-F_\theta (s)$.

\bp \label{prop: _geom_mean_between}
Let $p_\theta$, $\theta \in \Theta \subseteq \Rex $, be the probability mass function of a real observation $S$, having a continuous derivative with respect to $\theta$.
For each $s\in {\cal S}_0$, assume that  the function
\be \label{def_gamma_s}
\gamma_s (\theta) =  \left(\frac{\partial p_\theta(s)}{\partial \theta}\right)
/\left(-\frac{\partial F_\theta(s)}{\partial \theta}\right)=
\left(\frac{\partial p_\theta(s)}{\partial \theta}\right)\left/
h_s (\theta)\right.
\ee
 is decreasing on $\Theta$. Then $H_s(\theta) < H_s^G (\theta) < H_s^\ell(\theta)$ uniformly on $\Theta$.
\ep

The assumptions required in the previous proposition are satisfied by many important models. For example we have the following

\bc \label{cor:gamma_s_decreasing_NEF}
If $p_\theta$ is the probability mass function of a real NEF, then $H_s(\theta) < H_s^G (\theta) < H_s^\ell(\theta)$ uniformly on $\Theta$.
\ec

We now discuss the relationship between  $H_s^G$ and $H_s^A$.
\bp \label{prop:comparison_Hg_Ha}
Let $p_\theta$, $\theta \in \Theta \subseteq \Rex $, be the probability mass function of a real observation $S$, satisfying the following assumptions in addition to those stated in Proposition \ref{prop: _geom_mean_between}:
\be \label{properties_limit_gamma_s}
\lim_{\theta\rightarrow\inf \Theta} \gamma_s (\theta)= +\infty;
\hspace{5 mm}   \lim_{\theta\rightarrow\sup \Theta} \gamma_s(\theta)= -1.
\ee
Then, for each $s \in {\cal S}_0$, there exists $\theta^* \in \Theta$ (depending on $s$) such that   $H_s^G(\theta) <H_s^A (\theta)$ for $\theta < \theta^*$ and $H_s^G(\theta) > H_s^A (\theta)$ for $\theta \geq  \theta^*$.
\ep

The result in Proposition \ref{prop:comparison_Hg_Ha} is important in connection with confidence intervals, because it shows that $H_s^G$  gives, for a fixed  level, a  confidence interval smaller than that obtained from $H_s^A (\theta)$; see Figure 1 (graph 2) for an example.

Notice that the assumptions on $\gamma_s (\theta)$ in Proposition \ref{prop:comparison_Hg_Ha} are fulfilled by a real NEF with natural parameter space $\Theta=\Rex$, as it occurs in the binomial and Poisson models. However, these assumptions are not necessary to ensure the stated behavior of $H_s^G$ and $H_s^A$, that we conjecture to be quite general, as the following example shows.

\medskip
\noindent
\emph{Example 1}.
Consider an i.i.d. sample of size $n$ from a logarithmic distribution with parameter $\theta \in (0,1)$ with probability mass function
\[
p_\theta (x)=\frac{\theta^x}{-x\log(1-\theta)}I_{\{1,2,\ldots\}}(x).
\]
The  sufficient statistic $T=\sum_{i=1}^n X_i$ is distributed as
\[
p_\theta (t)=\frac{n! |s(t,n)| \theta^t}{t! (-\log(1-\theta))^n}I_{\{n,n+1,\ldots\}}(t),
\]
where  $s(t,n)$ is  the Stirling number of the first kind with arguments $t$ and $n$, see \citet{Johnson:2005}. The distribution of $T$ belongs to a real NEF with $F_\theta (t)$ decreasing in $\theta$, so that the fiducial distribution function $H_t$, for $t=n,n+1,\ldots$ and $\theta \in (0,1)$, is
\[
H_t(\theta)=1-F_\theta(t)=1-\sum_{j=n}^{t} \frac{n! |s(j,n)| \theta^j}{j! (-\log(1-\theta))^n}.
\]

\begin{figure}
\centering
\includegraphics[width=.49\textwidth,height=5cm]{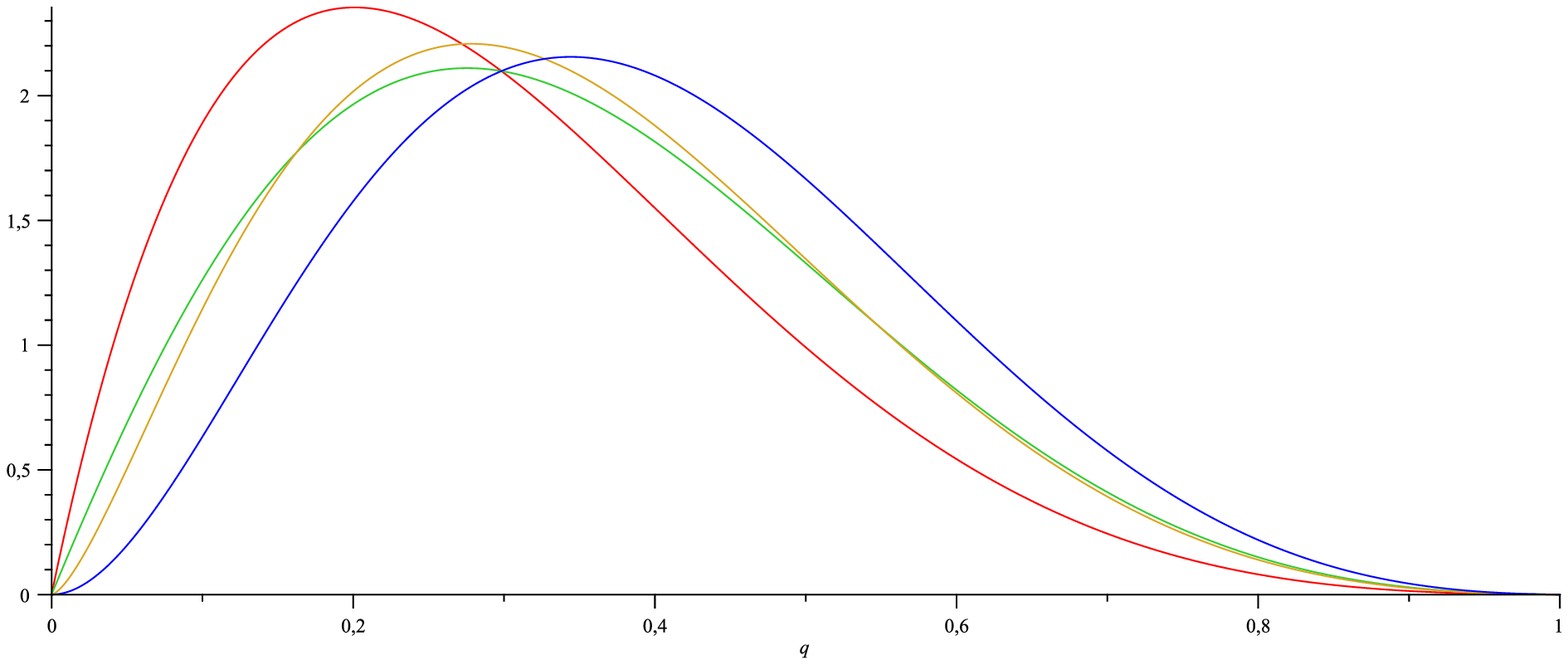}\hfil
\includegraphics[width=.49\textwidth,height=5cm]{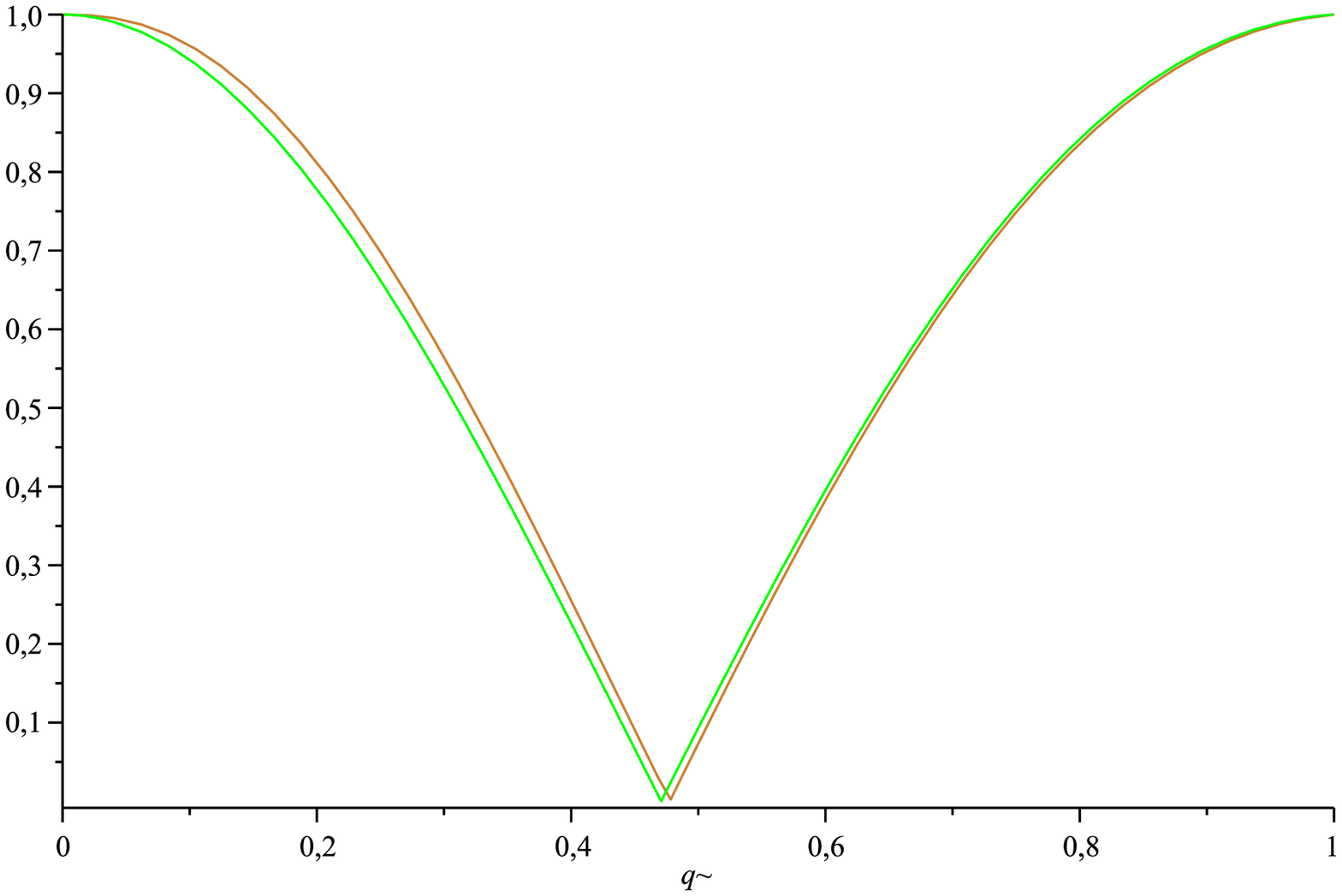}
\parbox{14cm}{
\caption{\small{Graph 1: Fiducial distributions for a sample from the logarithmic distribution ($n=10$, $t=12$): $h_t^{\ell}(\te)$ (red), $h_t^A(\te)$ (green), $h_t^G(\te)$ (yellow), $h_t(\te)$ (blue). Graph 2: Confidence curves for $h_t^A(\te)$ (green) and $h_t^G(\te)$ (yellow)).}}\label{fig:fid-log}}
\end{figure}

For this model
\[
\gamma_t(\theta) = -\left(\frac{\partial p_\theta(s)}{\partial \theta}\right)
/\left(\frac{\partial F_\theta(s)}{\partial \theta}\right)
= -\frac{|s(t,n)|\theta^{t-1}(n\theta+t(1-\theta)\log(1-\theta))}
{t!\sum_{j=n}^{t}|s(j,n)|\theta^{j-1}(n\theta+j(1-\theta)\log(1-\theta))/j!}.
\]
 It can be seen  that, for each $t \geq n$, $\gamma_t$ is decreasing in $\theta$ and
  \[
  \lim_{\theta \rightarrow 0^+}\gamma_t(\theta) = +\infty, \hspace{1cm}
   \lim_{\theta \rightarrow 1^-}\gamma_t(\theta) = -\left(\sum_{j=n}^t \frac{|s(j,n)|}{|s(t,n)|} \frac{t!}{j!} \right)^{-1}  \in (-1,0).
  \]
  Nevertheless, the fiducial distributions $H_t^G$ and $H_t^A$ behave as stated in Proposition  \ref{prop:comparison_Hg_Ha}, see Figure \ref{fig:fid-log} (graph 1).

Finally, we justify  our preference for $H^G_s$ versus $H^A_s$ showing that its confidence risk under quadratic penalty, as defined in \citet[Sec. 5.3]{Schweder:2016}, is uniformly better for all the important discrete models reported in Table \ref{tab_nef-qvf}. The confidence risk $R(\mu, H_s)$ for the mean parameter $\mu$ and a confidence (or  fiducial)  distribution $H_s$ under quadratic penalty is
\be
R(\mu,H_s)=\int{(\mu^\prime-\mu)^2 dH_s(\mu^\prime)}=E_\mu(\mbox{Var}^{H_s}(\mu))+E_\mu(\muh-\mu)^2,
\ee
where $\mbox{Var}^{H_s}(\mu)$ denotes the variance of $\mu$ under $H_s$, $E_\mu$ the expected value with respect to the distribution of $S$  and $\muh=E^{H_s}(\mu)$ is the mean of $\mu$ under $H_s$.
Now, recalling that for the binomial and the negative-binomial distribution in Table \ref{tab_nef-qvf}, assuming $m=1$ for simplicity,  we have $\mu=p$ and $\mu=(1-p)/p$,   it is easy to verify that for both these models and the Poisson model  $\muh$ is the same  under $H^G_s$ and $H^A_s$. As a consequence,  $R(\mu,H_s^A)-R(\mu,H_s^G)= E_\mu(\mbox{Var}^{H_s^A}(\mu))-E_\mu(\mbox{Var}^{H_s^G}(\mu))$ which becomes  $(4(n+1)(n+2))^{-1}$,  $(4(n-1)(n-2)^{-1}$ and $(4n^2)^{-1}$, for the three models above, respectively. All these values are  strictly positive for each $n$ uniformly in $\mu$.

\medskip

Let us now consider the fiducial distribution for a multivariate parameter defined in (\ref{fid-general}). For each discrete component of the product, starting from $\mbox{Pr}_{\f_{d-k+1}}\{T_k \leq t_k|\gT_{[k-1]}=\gt_{[k-1]},\gT_{-[d]} =\gt_{-[d]}; \gfi_{[d-k]}\}= F_{\f_{d-k+1}}(t_k|\gt_{[k-1]},\gt_{-[d]}; \gfi_{[d-k]})$ and
 $\mbox{Pr}_{\f_{d-k+1}}\{T_k < t_k|\gT_{[k-1]}=\gt_{[k-1]}, \gT_{-[d]}=\gt_{-[d]}; \gfi_{[d-k]}\}$, it is possible to  define a right and a left fiducial distribution, respectively, and hence  their geometric and arithmetic means. Notice that  each component of (\ref{fid-general}) involves a one-dimensional parameter and a real observation (the remaining quantities being fixed), so that the   Propositions \ref{prop: _geom_mean_between} and \ref{prop:comparison_Hg_Ha} can be applied.
 Multivariate fiducial distributions for discrete observations can thus be obtained combining in the various possible way these univariate distributions. In particular, we will consider
    $\Hgt(\gfi)$, obtained as the product of all the \emph{right} univariate conditional fiducial distributions,  $\Hgt^{\ell}(\gfi)$, obtained as the product of all the \emph{left} univariate conditional fiducial distributions,   $\Hgt^A(\gfi)$, defined as the product of the  $d$ mixtures $\Hgtk^A=(\Hgtk+\Hgtk^\ell)/2$ and finally $\Hgt^G(\gfi)$, corresponding to the density $h_{\boldsymbol{t}}^G(\gfi)$ obtained as the product of the $d$ geometric means $\hgtk^G\propto (\hgtk \cdot \hgtk^\ell)^{1/2}$. Notice that
$h_{\boldsymbol{t}}^G(\gfi)$  coincides with  the geometric mean of all the $2^{d}$   fiducial densities derived as described above.

\subsection{Fiducial inference and the sufficiency principle}\label{sec:generalized fiducial}

The step-by-step procedure introduced at the beginning of Section \ref{sec:step-by-step fiducial} gives a generalized fiducial distribution, according to  \citet{Hannig:2009}, if one considers as data-generating  equation $\gT=\gG(\gfi, \gU)$ with
\be
T_k=\left\{
\begin{array}{ll}
               G_k(\gfi,\gU_{[k]},\gU_{-[d]}) & k=1,\dots, d \\
               U_k & k=d+1,\dots, m
\end{array}
\right.,
\ee
where $\gU$ is a random vector with a completely known distribution.
The functions $G_k$ can be explicitly obtained iteratively as follows:
\be
T_1&=&G_1\left(\gfi,U_1,\gU_{-[d]}\right)=F^{-1}_{\f_{d}}\left(U_1|,\gU_{-[d]}; \gfi_{[d-1]} \right)\\
T_2&=&G_2\left(\gfi,U_1, U_2,\gU_{-[d]}\right)=F^{-1}_{\f_{d-1}}\left(U_2|G_1(\gfi,U_1,\gU_{-[d]}), \gU_{-[d]}; \gfi_{[d-2]}\right)
\ee
and so on.

It is interesting to observe that the generalized fiducial distribution $r(\te)$ given in \eqref{r-hannig} does not necessarily satisfy the sufficiency principle. This can be verified immediately looking at the Example 2 in \citet{Hannig:2013}, in which a uniform distribution on $(\te,\te^2)$ is considered and $r(\te)$ does not depend on the $X_i$'s only through the sufficient statistic $\gS=(X_{(1)},X_{(n)})$, where $X_{(i)}$ denotes the $i$-th order statistic.
Despite its simple form, this model is highly irregular, but
the inconsistency with the sufficiency principle of the generalized fiducial  distribution $r(\te)$   can also occur for more standard models. In particular, if  a real continuous sufficient statistic $S$ for a real parameter exists, one could derive two different fiducial distributions  starting from  $S$ or from the whole sample.  A simple example of this issue  can be easily constructed considering a  beta model with parameters 2 and $\theta$.
Another interesting example is the following.

\begin{figure}
\vspace{-1.5cm}
\centering
\includegraphics[width=.49\textwidth,height=5cm]{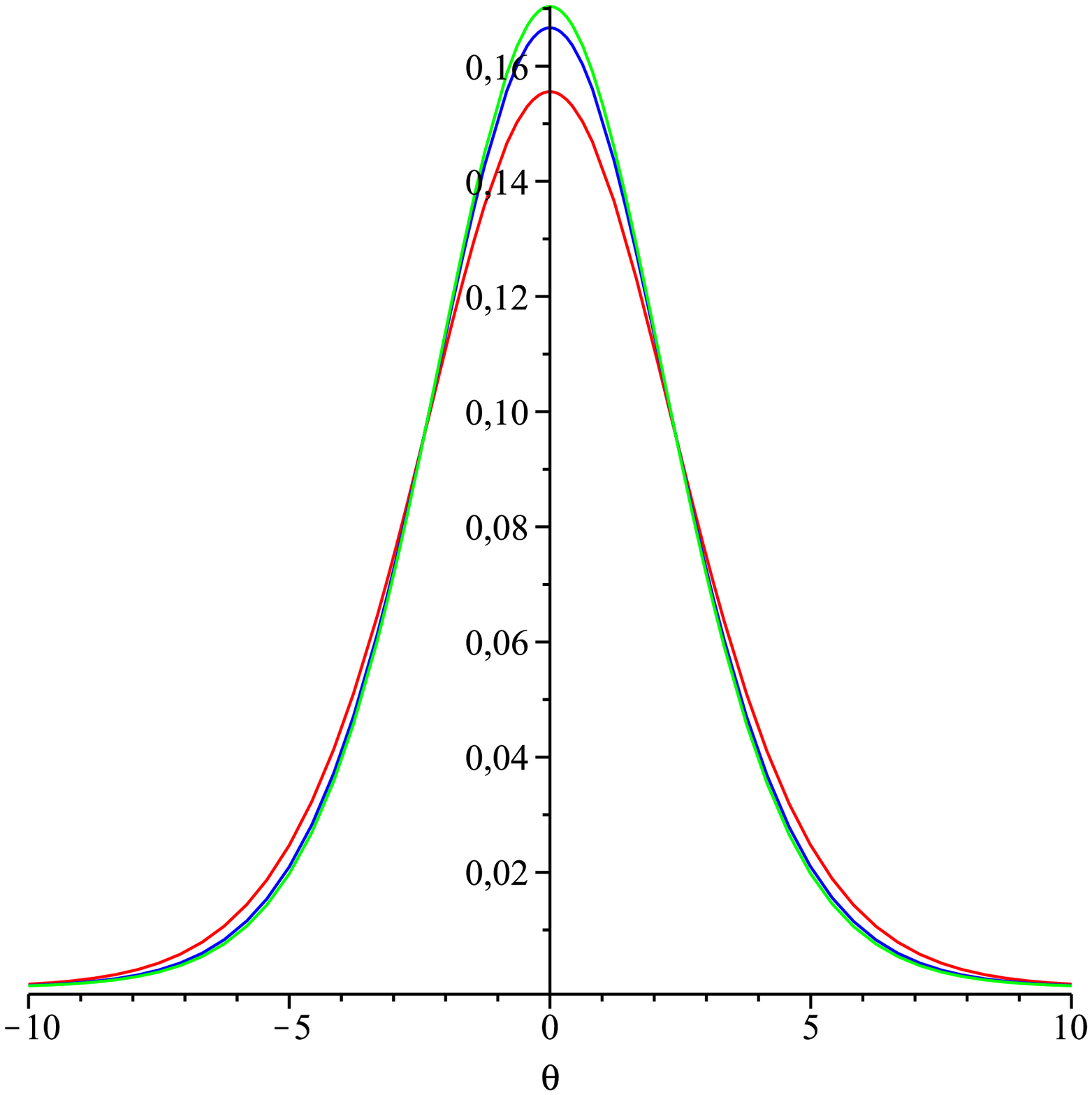}\hfil
\includegraphics[width=.49\textwidth,height=5cm]{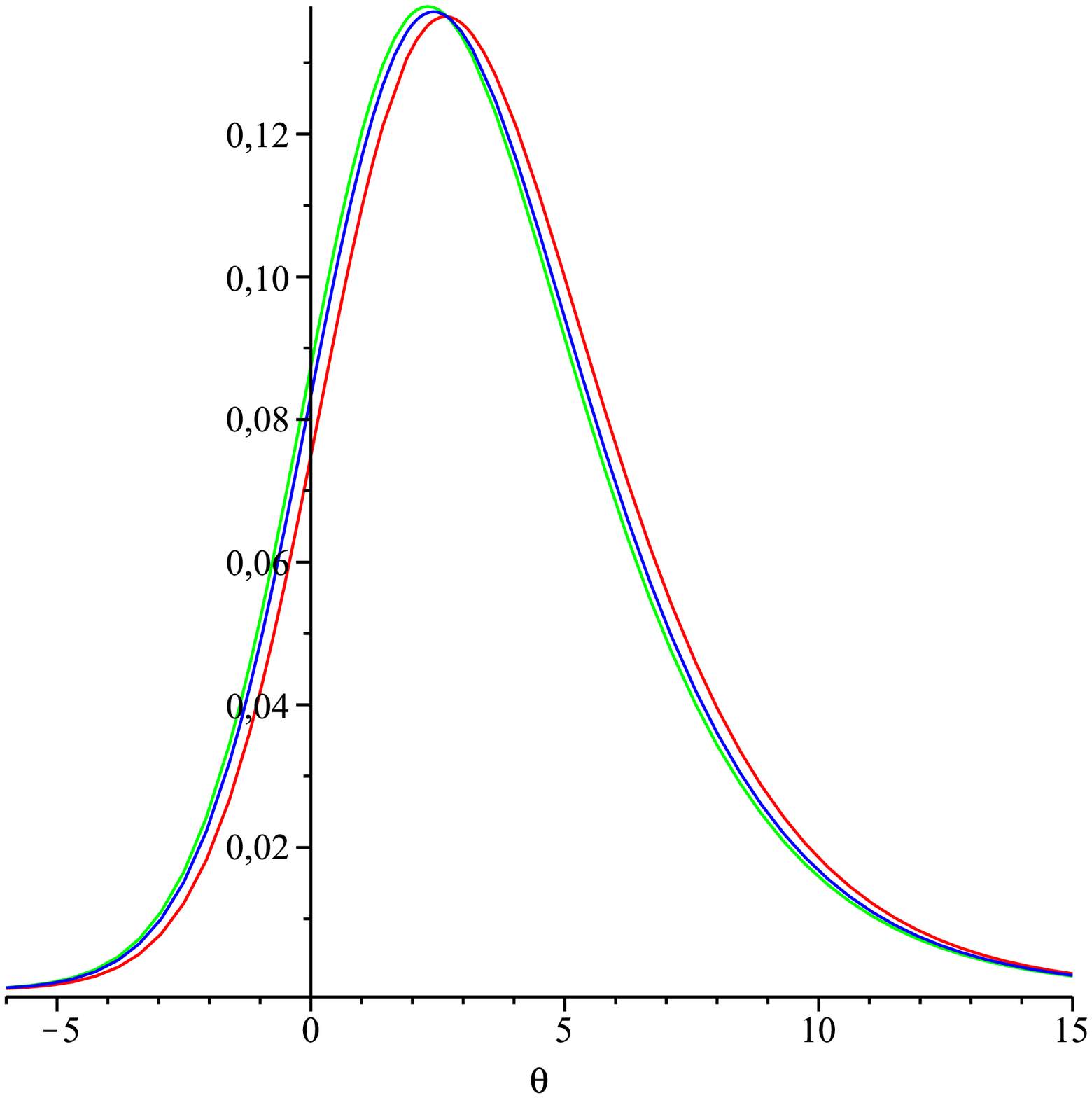}
\parbox{14cm}{
\caption{\small{Graph 1: Fiducial densities: $r(\te; x_1=0.05, x_2=0.95)$ (red);  $r(\te; x_1=0.5, x_2=0.5)$ (green); $h(\te; s=1)$ (blue). Graph 2: Fiducial densities $r(\te; x_1=0.02, x_2=0.48)$ (red);  $r(\te; x_1=0.2, x_2=0.3)$ (green); $h(\te; s=0.5)$ (blue).}}\label{fig:fid-exp-tronc}}
\end{figure}

\begin{figure}
\begin{center}
\includegraphics[width=9cm,height=4cm]{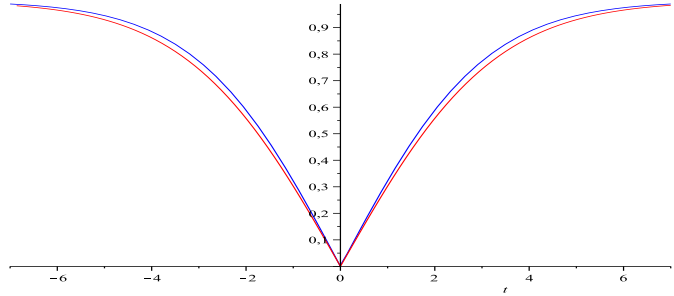}\\
\parbox{13cm}{
\caption{\small{Confidence curves for $r(\te; x_1=0.5,x_2=0.5)$ (red) and $h(\te; s=1)$ (blue). \label{fig:cc-exp-tronc}}}}
\end{center}
\end{figure}

\medskip
\noindent
\emph{Example 2}. Let $\gX=(X_1, \ldots, X_n)$ be an i.i.d. sample from a truncated  exponential density
$p_{\te}(x)=\te e^{-\te x}/(1-e^{-\te})$, $0<x<1$, $\te \in \Rex -\{0\}$. This density is not defined for $\te=0$, but it can be completed by continuity setting  $p_{0}(x)=1$.
The distribution function of $X_i$ is
$
F_{\te}(x_i)=(1-e^{-\te x_i})/(1-e^{-\te})$, $0<x_i<1$,
so that from \eqref{Hannig-J}  we have
\be
 J(\gx,\te)=\frac{s}{\te}+\frac{e^{-\te}}{\te(1-e^{-\te})} \sum_{i=1}^n(1-e^{-\te x_i}),
 \ee
 where $s=\sum_{i=1}^n x_i$. Thus, using \eqref{r-hannig}, we obtain
 \ben \label{r-exp-tronc}
 r(\te) \propto \frac{\te^{n-1}}{(1-e^{-\te})^{n+1}} e^{-\te s}  \left(s(1-e^{-\te})+e^{-\te} \sum_{i=1}^n (1-e^{-\te x_i})\right), \quad \te \in \Rex,
 \een
 which depends on the values of the specific $x_i$'s.
 Consider now the sufficient statistic $S=\sum_{i=1}^n X_i$ and, for simplicity, assume $n=2$. The density of $S$ is
 \be
 p_{\te}(s)=
 \left\{
\begin{array}{ll}
               \frac{\te^2}{(1-e^{-\te})^2}\; e^{-ts}s & \; 0<s\leq 1 \\
               \frac{\te^2}{(1-e^{-\te})^2}\; e^{-ts}(2-s) & \; 1<s<2
\end{array}
\right..
 \ee
and the generalized fiducial density \eqref{r-exp-tronc} reduces to $h_s(\te)=\partial F_{\te}(s)/\partial \te$.
In Figure  \ref{fig:fid-exp-tronc} we report the fiducial densities $r$ and $h_s$ for different values of $(x_1,x_2)$ and $s=x_1+x_2$. For $s=1$ all densities are symmetric with the mode in $0$, while  the dispersion is increasing in $|x_1-x_2|$, so that the more concentrated fiducial density is obtained for $x_1=x_2=0.5$. However, for $s\neq1$, the densities have different modes and are shifted to the left when $x_1$ increases. In all cases the fiducial density $h_s$ is in the middle of the various cases.
 Notice that $h_s$ has all the good properties discussed in \citet{Veronese:2014} and, in particular,  it is a confidence distribution because the model belongs to a NEF.
The confidence  intervals corresponding to $h_s(\theta)$ are slightly smaller than those corresponding to  $r(\theta)$, as can be seen from the confidence curves reported in Figure  \ref{fig:cc-exp-tronc}.
 For instance, when $x_1=x_2=0.5$, the 95\% confidence intervals are  (-4.191,4.191) and (-4.399,4.399), for $h_s(\theta)$ and $r(\theta)$,  respectively.

The  computation of  the fiducial distribution $\Hgt$ defined in \eqref{fid-general} is greatly simplified starting with the sufficient statistic instead of the whole sample.
However,  when  both the alternatives are feasible,  they seem  to lead  to the same result. In particular, the  following proposition states that the sufficiency principle is always satisfied by $\Hgt$ when there exists a complete sufficient statistic for the parameter.

\bp \label{prop:sufficiency}
Consider the fiducial distribution $\Hgt(\gfi)$, defined in \eqref{fid-general} and \eqref{fid-general-2}, with  $\gT$  a one-to-one transformation of the data $\gX=(X_1,\dots,X_n)$. If $\gS =\gS(\gT)$ is a complete and sufficient statistic of dimension $d$ for $\gfi$, such that
$\gS=\textbf{g}(\gT_{[d]}, \gT_{-[d]})$ is a one-to-one lower triangular transformation of $\gT_{[d]}$ for fixed $\gT_{-[d]}$,
 then the fiducial distribution $\Hgs(\gfi)$ for $\gfi$, obtained using $\gS$ instead of $\gT$ in \eqref{fid-general} and \eqref{fid-general-2}, coincides with $\Hgt(\gfi)$.
\ep

Notice that the   completeness of $\gS$ is not necessary to satisfy the sufficiency principle, as the following example shows.

\medskip

\noindent
\emph{Example 3}.
 Given an i.i.d. sample $\gX$ of size $n$ from a uniform distribution on $(\te,\te+1)$, it is immediate to verify that the sufficient statistic $\gS=(X_{(1)},X_{(n)})$ is not complete. Because $\te$ is a location parameter, $Z=X_{(n)}-X_{(1)}$ is an ancillary statistic and the fiducial distribution for $\te$ can be obtained starting from the distribution function of $X_{(n)}$ given $Z$, which is
$F_{\te}(x_{(n)}|z)= (x_{(n)}-z-\te)/(1-z)$,   $x_{(n)}-1<\te<x_{(n)}-z=x_{(1)}$.
Thus
\ben \label{fid-teta-teta+1}
h_{\gs}(\te) =- \frac{\partial}{\partial \theta} F_{\te}(x_{(n)}|z)= - \frac{\partial}{\partial \theta}  \frac{x_{(n)}-z-\te}{1-z}= \frac{1}{1-z}, \quad x_{(n)}-1<\te<x_{(n)}-z=x_{(1)}.
\een
If we start directly with $\gX$, we can consider the distribution function of $X_n$ given $\gZ=(Z_1,\dots,Z_{n-1})$, where $Z_i=X_n-X_i$.
Omitting tedious calculations, we have
\be
F_{\te}(x_n \mid \gz)=\frac{x_n-\te-\max(z_i,0)}{1+\min(z_i,0)-\max(z_i,0)}, \quad \te+\max(z_i,0)<x_n<\te+1+\min(z_i,0),
\ee
and thus, for $x_n-1-\min(z_i,0)<\te<x_n-\max(z_i,0)$,
\be
h_{\gx}(\te)=-\frac{\partial}{\partial \theta} F_{\te}(x_n|\gz)=\frac{1}{1+\min(z_i,0)-\max(z_i,0)}.
\ee
Observing that $\min(z_i,0)=z_i$ unless $x_n=x_{(n)}$ and recalling that $z_i=x_n-x_i$, $i=1,\dots, n-1$, it follows that $x_n-1-\min(z_i,0)=x_{(n)}-1$ and similarly $x_n-\max(z_i,0)=x_{(1)}$, so that $h_{\gx}(\te)$ coincides with $h_{\gs}(\te)$ given in \eqref{fid-teta-teta+1}.

\subsection{Conditionally reducible natural exponential families} \label{sec:cr-NEF}

Consider a multivariate natural exponential family whose density, with respect to a fixed $\sigma$-finite positive measure $\nu$, is given by
\ben \label{NEF}
p_{\gte}(\gx) =  \exp\left\{ \sum_{k=1}^d\te_k x_k  - M(\gte) \right \}, \;\; \gte=(\te_1\dots,\te_d) \in \gTe, \;\; \gx=(x_1,\dots,x_d) \in \Rex^d.
\een

A NEF is \emph{d-conditionally reducible}  (in the sequel cr-NEF) if its joint density  \eqref{NEF} can be factorized as a product of $d$ conditional densities each belonging to a real exponential family. More precisely if
  \ben \label{eqn:c_r_r}
   p_{\gfi}(\gx|\gfi(\gte))=
   \prod_{k=1}^{d}p_{\gfi_{k}}(x_k|\gx_{[k-1]}; \f_k(\gte))
   = \prod_{k=1}^{d}\exp\left\{\f_k(\gte) x_k-M_k(\f_k(\gte);\gx_{[k-1]})\right\},
  \een
 where $\gfi=(\f_1, \ldots,\f_d)$ is a one-to-one
function from $\gTe$ onto $\gfi(\gTe)=\gFi$. Furthermore, it can be shown
that $\gFi=\Fi_{1} \times \cdots \times \Fi_{d}$, with $\f_k \in
\Fi_{k}, k=1,\dots,d$, so that the $\f_{k}$'s are
variation independent.
Notice that  $\f_{k}$ is the natural
parameter of the $k$-th conditional distribution.
For details on these families, with emphasis on  \emph{enriched} conjugate priors and on reference Bayesian analysis, see \citet {Consonni:2001} and \citet{Consonni:2004}, respectively. Both these papers deal in particular with the families having \emph{simple quadratic variance function}, named NEF-SQVFs, which include, as most interesting cases, the multinomial and negative-multinomial models, see \citet{Casalis:1996} and Appendix A1.

\medskip

\noindent
\emph{Example 4} (Multinomial model).
Consider a  random vector $\gX$ distributed according to a multinomial distribution and denote by $p_k$ the probability of the $k$-th outcome $X_k$, $k=1,\dots,d$, with  $\sum_{k=1}^d x_k \leq N$, and $\sum_{k=1}^d p_k \leq 1$.  It is well know that the conditional distribution of $X_k$ given $\gX_{[k-1]}= \gx_{[k-1]}$, $k=2,
\ldots, d$, is Bi$(N-\sum_{j=1}^{k-1}x_j,
p_k/(1-\sum_{j=1}^{k-1}p_j))$, whereas the marginal distribution
of $X_1$ is Bi$(N,p_1)$. Since a binomial distribution is a real NEF, one can  factorize the multinomial distribution  as in (\ref{eqn:c_r_r}) with
\ben \label{eqn:binomial-fk}
 \phi_k=\log \frac{p_k}{1-\sum_{j=1}^{k}p_j}, \quad  M_k(\f_k; \gx_{[k-1]})=\left(N-\sum_{j=1}^{k-1}x_j\right) \log(1+e^{\f_k}), \quad \phi_k \in \Rex.
\een

\medskip

For models belonging to a cr-NEF, the construction of the fiducial distribution proposed in Section \ref{sec:step-by-step fiducial} drastically simplifies. The existence of  a sufficient statistic of the same dimension of the parameter makes the ancillary statistic not necessary, while the $\gfi$-parameterization, indexing each conditional distribution with a real parameter, implies the independence of the $\f_k$'s under the fiducial distribution.

\begin{prop}
\label{prop:fiducial-cr-NEF}
 Let $\gS$ be the sufficient statistic distributed according to a regular cr-NEF
on $\Rex^d$, parameterized by $\gfi \in \gFi=\Fi_{1} \times \cdots \times \Fi_{d}$, with $\Fi_{k}$ coinciding with the natural parameter space of the $k$-th conditional distribution.
Then, for $\gS=\gs$, with $s_k$, $k=1,\ldots, d$, satisfying  conditions similar  to those given before   \eqref{Hns},
\ben \label{fid-H-fi}
\Hgs(\gfi)=\prod_{k=1}^d \Hgsk(\f_k)=\prod_{k=1}^d (1-F_{\f_k}(s_k|\gs_{[k-1]})
\een
is a fiducial distribution function on $\gfi$ with density
\ben \label{fid-NEF-SQVF}
\hgs(\gfi)=\prod_{k=1}^d h_{\gs_{[k]}}(\f_k), \hspace{5 mm} \mbox{\emph{where}}   \hspace{3 mm}
h_{\gs_{[k]}}(\f_k)=\frac{\partial}{\partial \f} H_{\gs_{[k]}}(\f_k) = -\frac{\partial}{\partial \f_k} F_{\f_{k}}(s_k|\gs_{[k-1]}).
\een
\end{prop}

 The $\phi_k$'s are independent under $H_{\gs}(\phi)$ and thus their importance ordering is irrelevant. This fact also justifies the simplification in the index notation adopted in \eqref{fid-H-fi}. Notice, however, that the definition and the interpretation of the $\f_k$'s depend on the particular ordering considered for the $X_k$'s, as seen in  Example 4.

As recalled in Section \ref{sec:preli},  a general  definition of multi-dimensional confidence distribution does not exist. However, in our context, since  $H_{\gs}(\gfi)$ is constructed as a product of \emph{marginal} confidence distributions, it can be considered as a multivariate (possibly asymptotic) confidence distribution for $\gfi$.

 Some of the examples of Section \ref{sec:further examples} can be reconnected with this framework, but here we consider in specific the NEF-SQVFs, whose variance function is given in \eqref{basic}. For this class, with the exclusion of the negative-multinomial/hyperbolic secant distribution, it is possible to give a simple explicit expression of the fiducial density of $\gfi$, recalling the definition of $B_k(\f_k)$  given in \eqref{eqn:cond-cum}   and setting $z_k=z_{kk}$  in  \eqref{basic}. The specifications of $z_{kk}$ and $q$, appearing in \eqref{basic}, and of $B_k(\phi_k)$ can be found in Appendix A1.

\begin{prop} \label{prop:fid-NRF-SQVF}
Consider a sample of size $n$ from a Poisson/normal, a multinomial, or a negative-multinomial family on $\Rex^d$. If $\gS$ denotes the sufficient statistic, then the (right) fiducial distribution for $\gfi$ has
density
\ben \label{basic-fid}
  \hgs(\gfi) = \prod_{k=1}^d h_{\gs_{[k]}}(\f_k) \propto \prod_{k=1}^d \exp
  \left\{\f_k \left( s_k + z_k \right) - \left[n + q \left(\sum_{j=1}^{k-1} s_j-1 \right) \right] B_k(\phi_k)
  \right\},
 \een
while for the  negative-multinomial/gamma/normal family, with an $m$ dimensional negative multinomial component, the (right) fiducial distribution is given by
\ben \label{fiducial-NM/G/N}
h_{\gs}(\gfi) &=& \prod_{k=1}^d h_{[s_k]}(\f_k)  \propto
 \prod_{k=1}^{m}
  \exp \{\f_k (s_k + 1)\} (1-\exp(\f_{k}))^{n/q + \sum_{j=1}^{k-1} s_j-1}
  \nonumber \\
 & \times & \exp\left\{\f_{m+1}  s_{m+1} \right\}
(-\f_{m+1})^{n/q + \sum_{j=1}^{m} s_j-1 }
   \nonumber \\
 & \times & \prod_{k=m+2}^{d}  \exp\left\{\f_k s_k -
  n \: s_{m+1} \: \phi_k^2/2 \right\}.
\een
\end{prop}

Notice that the discrete components of a basic NEF-SQVF are integer-valued with $z_k=1$, so that the left fiducial distribution is obtained by the previous formulas replacing the term $(s_k+1)$ by $s_k$ in (\ref{basic-fid}) and (\ref{fiducial-NM/G/N}).
Thus it follows that the geometric mean $\hgs^G$ has the same structure  in  \eqref{basic-fid} and \eqref{fiducial-NM/G/N} with  $(s_k+1/2)$ instead of $(s_k+1)$.

\medskip

\noindent
{\em Example 4 (ctd.)}. For the multinomial family, because $q=-1/N$, $z_k=1$ and $B_k(\f_k)=N \log(1+e^{\f_k})$, $k=1,\ldots,d$, it easily follows from  formula \eqref{basic-fid}   that
 \ben \label{multinomial-Gfid-fi}
  \hgs^G(\gfi)  = \prod_{k=1}^d h^G_{\gs_{[k]}}(\f_k)
  = \prod_{k=1}^d
  \left\{\frac{1}{\texttt{B}(s_k+\frac{1}{2},nN-\sum_{j=1}^{k} s_j+\frac{1}{2})}\frac{e^{\f_k (s_k + \frac{1}{2})}}{\;(1+e^{\f_k})^{nN - \sum_{j=1}^{k-1} s_j+1} } \right\},
 \een
 where $\texttt{B}(\cdot,\cdot)$ denotes the beta function.

\medskip

The fiducial distribution for $\gfi$ not always is of particular interest in itself, but  it can be used as a starting point for the construction of the fiducial distribution for alternative and more relevant parameters.
%
We consider here the mean-parameter $\gmu$, which is a lower triangular transformation of $\gfi$, see \eqref{eqn:mu-phi}, so that its fiducial distribution can be directly obtained from that of $\gfi$ thanks to Proposition \ref{prop:invariance}.
  \begin{cor} \label{cor:fid-mu-sqvf}
  The (right) fiducial distribution for the mean parameter $\gmu$, relative to the ordering $\mu_1, \ldots,\mu_d$, for the following NEF-SQVFs on $\Rex^d$, has density:
  \begin{itemize}
\item Poisson/normal family (with $m$ Poisson components)
   \ben \label{fid-poisson-mu}
   h_{\gs}(\gmu)\propto \prod_{k=1}^m \mu_k^{s_{k-1}}\exp(-n\mu_k)  \prod_{k=m+1}^d \exp\left\{-\frac{n}{2 \s^2}(\mu^2-2\mu_k s_k/n)\right\},
   \een
   which corresponds to the product of $m$   densities  Ga$(s_k,n)$, $k=1,\ldots,m$ and  $(d-m)$ densities  N$(s_k/n, \s^2/n)$, $k=m+1,\ldots,d$.
\item Multinomial family
\ben \label{fid-multinomial mu}
h_{\gs}(\gmu) \propto \prod_{k=1}^d \mu_k^{s_k} \left(N-\sum_{j=1}^k \mu_j\right)^{\gamma_k},
\een
where $\gamma_k=-1$ for $k=1,\dots,d-1$ and $\gamma_d=Nn-1-\sum_{j=1}^d s_j$.
\item Negative-multinomial family, with $R$ occurrences in the $(d+1)$-th cell
\ben \label{fid-neg-multinomial-mu}
h_{\gs}(\gmu) \propto \prod_{k=1}^d \mu_k^{s_k} \left(R+\sum_{j=1}^k \mu_j\right)^{\gamma_k},
\een
where $\gamma_k=-1$ for $k=1,\dots,d-1$ and $\gamma_d=-Rn-1-\sum_{j=1}^d s_j$.
\item Negative-multinomial/gamma/normal family  (with an $m$-dimensional negative - multinomial component with $R$ occurrences in the $(m+1)$-th cell)
\ben \label{fid-neg-multinomial-gamma-normal-mu}
h_{\gs}(\gmu) &\propto& \prod_{k=1}^m \mu_k^{s_k} \left(R+\sum_{j=1}^k \mu_j\right)^{\gamma_k}\\
&\times& \left(R+\sum_{j=1}^m \mu_j\right)^{Rn+\sum_{j=1}^m s_j}        \mu_{m+1}^{-(Rn+\sum_{j=1}^m s_j)-1} \exp\left\{-\frac{s_{m+1}\left(R+\sum_{j=1}^m \mu_j\right)}{\mu_{m+1}}\right\} \nonumber\\
&\times& \prod_{k=m+2}^d \mu_{m+1}^{-(d-m-1)}\exp\left\{-\frac{n s_{m+1}}{2 \mu^2_{m+1}}\left(\mu_k^2-2\,\mu_k \mu_{m+1}\frac{s_k}{n}\right)\right\}, \nonumber
\een
where $\gamma_k=-1$ for $k=1,\dots,m-1$ and $\gamma_m=-Rn-1-\sum_{j=1}^m s_j$.
Notice that the density of $\mu_{m+1}$ given $\gmu_{[m]}$ is an In-Ga$(Rn+\sum_{j=1}^m s_j, s_{m+1}(R+\sum_{j=1}^m \mu_j))$, while the density of  $\mu_k$ given $\gmu_{[k-1]}$, $k=m+2, \dots,d$, is a N$(\mu_{m+1}s_k/n, \mu_{m+1}^2/(n s_{m+1}))$ depending only on $\mu_{m+1}$.
\end{itemize}
\end{cor}

\noindent
\emph{Example 4 (ctd.)}. Inference for the multinomial distribution  is usually performed for the cell-probabilities parameter $\gp=(p_1\,\dots,p_d)$.  Since $p_k=\mu_k/N$, the  fiducial distribution $h^G$ for $\gp$ is easily derived from \eqref{fid-multinomial mu}, noting that
the \emph{left} fiducial density can be obtained replacing  $(s_k+1)$ by  $s_k$ in \eqref{basic-fid}, (and not in \eqref{fid-multinomial mu}, which is derived aggregating the hyperparameters).
It follows that the geometric mean $h_{\gs}^G(\gp)$ is given by
 \ben \label{multinomial-fid-p}
h_{\gs}^G(\gp) \propto \prod_{k=1}^d p_k^{s_k-1/2} \left(1-\sum_{j=1}^k p_j\right)^{\gamma_k}, \quad \quad \sum_{k=1}^d p_k=1, \quad 0<p_k<1,
\een
with $\gamma_k=-1/2$ for $k=1,\dots,d-1$ and $\gamma_d=Nn-1/2-\sum_{j=1}^d s_j$. This is a generalized Dirichlet  distribution. Clearly, $h_{\gs}^G(\gp)$ in \eqref{multinomial-fid-p} refers to the specific order of importance  $p_1, p_2,\ldots,p_d$. If we change this order, the fiducial distribution will change accordingly.

\medskip

Similarly, for the negative-multinomial model, with $R$ occurrences in the $(d+1)$-th cell, $h_\gs^G (\gfi)$ can be easily computed from (\ref{basic-fid}) observing that $z_k=1$, $q=1/R$ and $B_k(\f_k)=-R \log (1-\exp(\f_k))$.

\section{Connections with objective Bayesian inference} \label{sec:connection reference}

As mentioned in Section \ref{sec:_intro}, if we look at  fiducial inference as a way to obtain a  distribution on  the parameter space of the model without any prior information, it appears natural to compare it with objective Bayesian inference.  Recall that when a fiducial distribution coincides with a posterior, the corresponding prior is called  \emph{fiducial prior}.

The step-by-step construction of the fiducial distribution $h_{\gt}(\gfi)$ defined  in  \eqref{fid-general} is based on the inferential importance ordering of the parameter components $\f_1, \ldots, \f_d$.  This aspect is also  crucial in the procedure adopted to construct reference priors,  see
\citet[Sec. 5.4.5]{Bernardo:1994}. The reference prior $\pi^R$ for a parameter  $\gfi$ is generated by successive conditioning,  established by the importance ordering of its components,  as $\pi^R(\gfi)=\prod_{k=1}^d \pi^R(\f_{d-k+1}|\gfi_{[d-k]})$.  It is widely recognized  that the dependence of the reference prior on the choice of the parameter of interest is necessary to obtain good frequentist properties such as coverage  and consistency. For a one-dimensional parameter $\f$   the reference prior coincides with the Jeffreys prior $\pi^J (\f) \propto I(\f)^{1/2}$, where $I(\f)$ denotes the Fisher information.
While the Jeffreys prior is invariant under a reparameterization of the model, the reference prior (and thus the reference posterior)  is  generally not invariant unless the transformation from $\gfi$ to $\gla=(\la_1, \ldots, \la_d)$ is  lower triangular,  see \citet{Datta:1996}.
Thus the reference posterior has the same invariance property of the fiducial distribution proved in Proposition \ref{prop:invariance}.

Recently \citet{Berger:2015} recognize the existence of situations in which one is interested simultaneously in all the parameter components of the model, or in none of them  but a prior (and thus a posterior) distribution is necessary to perform other inferences such as predictions. In these cases an \emph{overall} prior is needed. Its determination is an open problem but, as they  highlight, when there exists a \lq\lq common reference prior for all parameters\rq \rq{}, this is the natural choice for the overall prior. A similar problem occurs in our context and we will comment on this aspect in the following sections.
Notice that here the fiducial distribution \eqref{r-hannig} suggested by Hannig can be a good choice.

\subsection{Location-scale parameter models} \label{sec:location-scale}

For location-scale parameter models the fiducial prior exists and coincides with the reference prior.
Assume  first that only one parameter, $\te$, is unknown. In this case the model admits an ancillary statistic $\textbf{Z}$ and, in particular, we  take  $Z_i=X_i-X_1$ or $Z_i=X_i/X_1$, $i=2,\dots, n$, if $\te$ is a location or a scale parameter, respectively.

\begin{prop} \label{prop:fid location or scale}
Let $\gX=(\vecX)$ be an i.i.d. sample  from a density $p_\te$, $\te \in \Te\subseteq \Rex$.
If $\te$ is a location or a scale parameter, then the fiducial distribution coincides with the Bayesian posterior obtained with the Jeffreys prior $\pi^J(\te) \propto 1$ or $\pi^J(\te) \propto 1/\te$, respectively.
 \end{prop}

\noindent
\emph{Example 5.} Let $\gX$ be an i.i.d. sample from the  uniform distribution on $(0,\te)$, $\te>0$, so that $\theta$ is a scale parameter.
 First notice that $S=X_{(n)}$ is a sufficient statistic for $\theta$ and thus we can obtain directly the fiducial distribution
 \ben \label{F-max-uniform}
 h_s(\te)=\frac{\partial}{\partial \te} H_s(\te)=-\frac{\partial}{\partial \te} F_\te(s)=-\frac{\partial}{\partial \te} \left(\frac{s}{\te}\right)^n = \frac{n s^n}{\te^{n+1}} , \quad \te > s.
 \een
However the same result can be obtained without resorting to the sufficient statistic. Set $w=\mbox{max}(z_2, \ldots, z_n)$ and consider the distribution function of $X_1$ given the ancillary statistic $\gZ = (X_2/X_1, \ldots,$ $ X_{n}/X_1)$
\ben \label{fid_uniforme}
F_\te(x_1|\gz)= \left\{ \begin{array}{ll} \left(\frac{x_1}{\te}\right)^n & \hspace{5mm} 0<x_1<\te, \quad 0<w \leq 1\\
\left(\frac{x_1 w}{\te}\right)^n & \hspace{5mm} 0<x_1<\te, \quad w>1
\end{array} \right..
\een
Now, because  $w\leq1$ means $x_1=\mbox{max}(\vecx)$, while for $w>1$ we have $x_1w=\mbox{max}(x_2, \ldots, x_n)$, expression \eqref{fid_uniforme}, as a function of $\theta$, is equivalent to $F_\te(s)$ appearing in \eqref{F-max-uniform} and thus provides the same fiducial distribution. It is immediate to verify that it coincides with the Jeffreys posterior.

\medskip

A case in which the sufficient statistic is not one-dimensional and thus it is necessary to use an ancillary statistic can be found  in the previous Example 3. Trivially $h_{s}(\te)$ given in  \eqref{fid-teta-teta+1} coincides with the Bayesian posterior obtained by $\pi^J(\te) \propto 1$.

Consider now a model with a location parameter  $\te$ and a scale parameter $\s$, both unknown. Given an i.i.d. sample of size $n$, an ancillary statistic  is, for example, $\gZ=(Z_3,\ldots,Z_n)$, with $Z_j=(X_j-X_1)/Z_2$, $j=3, \ldots,n$, where $Z_2=X_2-X_1$ is marginally ancillary for $\te$. Then, the one-to-one transformation from  $\gX$ to $(X_1,Z_2,\gZ)$  allows to write
the sampling distribution as $ p_{\s}(z_2|\gz) p_{\te}(x_1|z_2,\gz; \s)p(\gz)$.
Note that in specific contexts other  transformations could be more appropriate. For example, in a normal model one could  use $(\bar{X}=\sum_{i=1}^n X_i/n$, $S^2=\sum_{i=1}^n (X_i-\bar{X})^2, \gZ)$ with $Z_j=(X_j-\bar{X})/S$, $j=3,\dots,n $, so that
the factorization becomes $p_{\s}(s^2|\gz) p_{\te}(\bar{x}|s^2,\gz;\s) p(\gz)$.

\begin{prop} \label{prop:fid location-scale}
Let $\gX=(\vecX)$ be an i.i.d. sample from a density $p_{\te,\s}$, where $\te$ and $\s$ are  a location and a scale parameter, respectively.
Then the fiducial distribution  $h_{\gx}(\sigma,\theta)$  for $(\s,\te)$ coincides with the Bayesian posterior obtained with the reference prior $\pi^R_{\s,\te}(\s,\te ) \propto 1/\s$.
 \end{prop}

Notice that $\pi^R(\s,\te) \propto 1/\s$ is different from  $\pi^J(\s,\te) \propto 1/\s^{2}$ obtained by the Jeffreys rule which, as already recalled, is not suitable for multidimensional parameters. Furthermore, while $\pi^R$  does not depend on the ordering of $\te$ and  $\s$, the step-by-step fiducial distribution is in general not allowable if  the ordering is reversed. However, $h_{\gx}(\s,\te)$ coincides with the fiducial distribution  obtained through other \lq \lq symmetric\rq \rq{} approaches, see \citet{Hannig:2009} and \citet{Fraser:1961}. Thus  the inferential ordering  of importance seems irrelevant for this model and $ h_{\gx}(\s,\te)$ can be assumed as an overall fiducial distribution.


\subsection{Exponential families} \label{sec:bayesian-NEFs}

 Lindley (1958) was the first to study the existence of a fiducial prior, analyzing in particular the case of  continuous real NEFs and  proving that it exists only for gaussian (with known variance) and gamma (with known shape) models.
A full characterization of  the  real NEFs  which admit  a fiducial prior is given in \citet{Veronese:2014}. The following proposition summarizes their results.
\begin{prop}
\label{prop:FID_JEF}
 Let $\cal F$ be a real NEF with natural parameter $\theta$.
 \begin{itemize}
   \item[i)] A fiducial prior exists if and only if $\cal F$ is  an affine transformation of one of the following families: normal with known variance, gamma with known shape parameter, binomial, Poisson and negative-binomial. For the three discrete families, the fiducial prior exists for all $H_s$, $H_s^\ell$ and  $H_s^G$.
   \item[ii)] When a fiducial prior exists, it  belongs to the family of conjugate distributions. Moreover, it coincides with the Jeffreys prior for   continuous NEFs and for discrete  NEFs too if we choose $H_s^G$  as the fiducial distribution.
   \item[iii)] The fiducial distribution  $H_s$ (or $H_s^A$ in the discrete case)  and the Bayesian posterior distribution corresponding to the Jeffreys prior have the same Edgeworth's expansion  up to the term of order $n^{-1}$.
 \end{itemize}

\end{prop}
The previous results establish a strong connections between  Jeffreys posteriors and fiducial distributions for real NEFs, and thus the two different approaches lead, in some sense, to the same \emph{objective inference}. A discussion about the coverage of the fiducial and the  Jeffreys intervals and their good frequentist properties, in particular when compared with the standard Wald intervals, is given in \citet[Section 5]{Veronese:2014}.

Consider now a cr-NEF.
 It is easy to verify that the fiducial distribution $h_{\gs}(\gfi)$  in  \eqref{basic-fid}  belongs to the enriched conjugate family defined in \citet[Section 4.3]{Consonni:2001}. This fact is the key-point to prove the following proposition.

\begin{prop} \label{prop: fid-prior-cr-NEF}
Let $\gS$ be a sufficient statistic  distributed according to a cr-NEF on $\Rex^d$, parameterized by $\gfi=(\f_1, \ldots,\f_d)$. Then   a fiducial prior for $\gfi$ exists  if and only if the conditional distribution of $S_k$ given $\gS_{[k-1]}=\gs_{[k-1]}$ is an affine transformation of one of the following families: normal with known variance, gamma with known shape parameter, binomial, Poisson and negative-binomial.

In particular, all basic NEF-SQVFs, with the exclusion of the Negative-multinomial/ hyperbolic secant, admit a fiducial prior, which belongs to the enriched conjugate family. Moreover, if for the discrete components of these models we consider the geometric mean $h_{\gs_{[k]}}^G$, then the product of the Jeffreys priors computed from the conditional distribution of $S_k$ given $\gS_{[k-1]}=\gs_{[k-1]}$,  the reference prior and the fiducial prior are all equal.
\end{prop}

\noindent
{\em Example 4
 (ctd.)}. The multinomial distribution is a basic NEF-SQVF and thus from Proposition \ref{prop: fid-prior-cr-NEF}, setting $s_k=n=0$ in $\hgs^G(\gfi)$ given in  \eqref{multinomial-Gfid-fi}, we obtain  the fiducial prior
 \ben \label{fid-prior fi}
 \pi(\gfi) \propto \prod_{k=1}^d e^{\f_k/2}/(1+e^{\f_k/2}),
  \een
  which coincides with the reference prior and with the product of the Jeffreys priors for $\f_{k}$, $k=1,\ldots,d$, computed on the distribution of $X_k$ given $\gX_{[k-1]}=\gx_{[k-1]}$.
\medskip

Finally,  we observe that the fiducial distribution  \eqref{fid-NEF-SQVF}  is always an overall fiducial distribution for $\gfi$. However,  the $\gfi$-parameterization is often not interesting in itself even if in some cases it is strictly related with a more relevant one. For example,
following \citet{Berger:2015}, consider a multinomial model applied to  directional data, as it happens for outcomes from an attitude survey. In this case the cells are naturally ordered, so that it is meaningful to reparameterize the model in terms of the conditional probabilities $p^*_k=\exp(\phi_k)/(1+\exp(\phi_k))$, $k=1, \dots,d$. Then $\pi(\gfi)$ in \eqref{fid-prior fi} induces on  $\gp{^*}=(p^*_1, \dots,p^*_d)$ an overall fiducial prior which is a product of independent Be(1/2,1/2) distributions coinciding with the overall reference prior.

\section{Further examples} \label{sec:further examples}

\subsection{Examples concerning normal models}
i) \emph{Difference of means.} Consider two independent normal i.i.d. samples, each of size $n$, with known common variance $\s^2$ and means $\mu_1$ and $\mu_2$, respectively. The sufficient statistics  are the sample sums  $S_1$ and $S_2$, with $S_i \sim$  N($n\mu_i, n\s^2$), $i=1,2$.
 If the parameter of interest is  $\phi_1=\mu_2-\mu_1$,
 we  can reparameterize the joint density of $(S_1, S_2)$ in $(\phi_1=\mu_2-\mu_1, \phi_2=\mu_1)$,
so that
%
the conditional distribution of $S_2$ given $S_1+S_2$, being N($(n\phi_1 + s_1+s_2)/2, n\s^2/2$), depends only on $\phi_1$.  From Table \ref{tab_nef-qvf},  the fiducial distribution of $\phi_1/2+(s_1+s_2)/(2n)$ is N$(s_2/n, \s^2/(2n))$, and thus  $\phi_1$ is N$(\bar{x}_2-\bar{x}_1, 2\s^2/n)$, where $\bar{x}_i=s_i/n$.
Because $S_1+S_2$ is N$(\phi_1+2\phi_2, 2n\s)$, arguing as before,
the fiducial distribution of $\phi_2$ given $\phi_1$ is  N($(\bar{x}_1+\bar{x}_2 - \phi_1)/2, \s^2/(2n)$), so that $h_{S_1,S_2}(\phi_1,\phi_2)=h_{S_1,S_2}(\phi_1)h_{S_1+S_2}(\phi_2|\phi_1)$. Notice that the same joint fiducial distribution is obtained if we consider the ordering $(\phi_2, \phi_1)$ or even if we compute the (marginal) fiducial distributions of $\mu_1$ and $\mu_2$ and obtain that of $(\phi_1,\phi_2)$ through the change-of-variable rule. Thus the ordering of the parameter is irrelevant and   $h_{S_1,S_2}(\phi_1,\phi_2)$ is an overall fiducial distribution. Furthermore, it coincides with the reference posterior obtained with a constant prior and  the marginal distribution of $\phi_1$ and $\phi_2$  are both confidence distributions.

\medskip

ii) \emph{Many normal means (Neyman Scott Problem).}
Consider $n$  samples of size two $(X_{i1}, X_{i2})$, with each $X_{ij}$ independently distributed according to a N$(\mu_i, \s^2)$, $i=1,\dots,n$ and let $\bar{X}_i=(X_{i1}+X_{i2})/2$ and $W=\sum_{i=1}^n(X_{i1}-X_{i2})^2$.
The aim is to make inference on the common variance $\s^2$, with nuisance parameter $\gmu=(\mu_1,\ldots, \mu_n)$. This well known example  is used to show that the   maximum likelihood estimator  $\hat{\sigma}^2=W/(4n)$ of $\sigma^2$ is inconsistent, because  $W/(4n) \rightarrow \s^2/2$,  $n \rightarrow \infty$. To obtain the fiducial distribution of $\s^2$, first notice that  the joint distribution of the sufficient statistics $\bar{\gX}=(\bar{X}_1, \ldots, \bar{X}_n)$ and $W$ can be factorized as $ \left(\prod_{i=1}^n p_{\mu_i,\s^2}(\bar{x_i})\right)p_{\s^2}(w)$,
for the independence of  $\bar{\gX}$ and $W$, with $W\sim$ Ga$(n/2,1/(4\s^2))$.
Using  Table \ref{tab_nef-qvf} one can easily obtain from $p_{\s^2}(w)$ the fiducial distribution for $1/(4\s^2)$, and hence that for $\s^2$ which is In-Ga$(n/2,w/4)$, while that of each $\mu_i$ given $\s^2$, derived from $p_{\mu_i,\s^2}(\bar{x_i})$,  is   N$(\bar{x_i}, \s^2/2)$. As a consequence
$h_{\bar{\gx},w}(\s^2,\gmu)=\left(\prod_{i=1}^n h_{\bar{x}_i}(\mu_i|\s^2)\right) h_{w}(\s^2)$.
This distribution coincides with the posterior obtained from the order invariant reference prior  $\pi^R(\s^2, \mu_1, \ldots, \mu_n)$ $\propto 1/\s^2$ and does not present the inconsistency of the likelihood estimator, which instead occurs for the posterior distribution obtained from the Jeffreys prior $\pi^J(\s^2, \mu_1, \ldots, \mu_n)\propto 1/\s^{n+2}$.

\subsection{Comparison of two Poisson rates} \label{sec: comparison Poisson}

The comparison of Poisson rates  $\mu_1$ and $\mu_2$ is a classical  problem arising in many contexts, see for example \citet{Lehmann:2005} for a discussion on an unbiased uniformly most powerful test for the ratio $\f_1=\mu_2/\mu_1$.
Given two i.i.d. samples of size $n$ from two independent Poisson distributions, the sufficient statistics are the sample sums $S_1$ and $S_2$,  with $S_i \sim$ Po($n\mu_i$), $i=1,2$. Reparameterizing the joint density of $(S_1, S_2)$ in $(\f_1=\mu_2/\mu_1, \f_2=\mu_1+\mu_2)$, we have that
the conditional distribution of $S_2$ given $S_1+S_2$ is Bi($s_1+s_2,  \f_1/(1+\f_1)$) and the marginal distribution of $S_1+S_2$ is Po($n\f_2)$. Thus the sampling distribution is a cr-NEF and we can apply \eqref{fid-H-fi}. Using Table  \ref{tab_nef-qvf},  the fiducial density for $\f_1/(1+\f_1)$, derived from the conditional distribution of  $S_2$ given $S_1+S_2$  is Be$(s_2+1/2, s_1+1/2)$ which implies
\ben \label{eqn:fid-ratio-poisson}
h^G_{s_1,s_2}(\phi_1)= \frac{1}{\texttt{B}(s_2+1/2,s_1+1/2)}\phi_1^{s_2-1/2}(1+\phi_1)^{-s_1-s_2-1}, \hspace{5 mm} \phi_1>0.
\een
From the marginal distribution of $S_1+S_2$ and using again Table  \ref{tab_nef-qvf},  it follows that $h_{s_1+s_2}^G(\f_2)$ is  Ga($s_1+s_2+1/2, n$) and thus $h_{s_1,s_2}^G(\f_1,\f_2)=h_{s_1,s_2}^G(\f_1)h_{s_1+s_2}^G(\f_2)$.
This joint  fiducial distribution is order-invariant, coincides with the reference posterior according to Proposition \ref{prop: fid-prior-cr-NEF}, and is an overall distribution for $(\f_1, \f_2)$. Notice that  $h_{s_1,s_2}^G(\f_1)$ is a confidence distribution and that it  differs from the fiducial distribution induced on $\f_1$ by the two independent marginal fiducial densities for $\mu_1$ and  $\mu_2$.

\subsection{Bivariate binomial}

A Bayesian  analysis for the bivariate binomial model has been discussed by \citet{Crowder:1989} in connection with a microbiological application.
Consider $m$ spores, each with a probability $p$ to germinate, and denote by $R$ the random number of germinating spores, so that $R$ is Bi($m,p$). If $q$ is the probability that one of the latter spores bends in a particular direction and $S$ is the random number of them, the probability distribution of $S$ given $R=r$ is Bi($r,q$). The joint distribution of $R$ and $S$ is called \emph{bivariate binomial}.
Crowder and Sweeting  observe that the Jeffreys prior $\pi^J(p,q) \propto p^{-1} (1-p)^{-1/2}q^{-1/2} (1-q)^{-1/2}$ is not satisfactory for its asymmetry in $p$ and $1-p$, while
\citet{Polson:1990} show that this fact does not occur using the order-invariant reference prior  $\pi^R(p,q) \propto p^{-1/2} (1-p)^{-1/2}q^{-1/2} (1-q)^{-1/2}$ which is the product of the two independent Jeffreys priors.

The joint fiducial density $h^G_{r,s}(q,p)$ can be obtained as the product of $h^G_{r,s}(q)$, derived from  the conditional model Bi($r,q$) of $S$ given $R=r$, and $h^G_{r}(p|q)$, derived from the marginal model Bi($m,p$) of $R$, which does not depend on $q$. Thus $p$ and $q$ are independent under  $h^G_{r,s}$ so that it is an overall  fiducial distribution. Because for the binomial model  the fiducial prior  is equal to the Jeffreys prior, see Proposition \ref{prop:FID_JEF}, it follows immediately that $h^G_{r,s}(q,p)$ coincides with the reference posterior.
All previous conclusions hold even if we consider the alternative  parametrization  $(\eta=pq, \la=p(1 - q)/(1 -pq))$.

\subsection{Ratio of parameters of a trinomial distribution}

 \citet{Bernardo:1998}  perform the Bayesian reference analysis for the ratio of two multinomial parameters presenting some applications. In particular they discuss the case of  $(X_1,X_2)$ distributed according to a trinomial distribution with parameters $n$ and $\gp=(p_1,p_2)$, and provide the joint reference prior for $(\phi_1=p_1/p_2, \phi_2=p_2)$,  with  $\phi_1$ the parameter of interest. Then they derive the marginal reference posterior for $\f_1$ which is
\ben \label{eqn:ref_eta}
\pi^R(\phi_1|x_1,x_2) \propto \phi_1^{x_1-1/2}(1+\phi_1)^{-x_1-x_2-1}.
\een
To find the fiducial distribution of $\phi_1$, we reparameterize the trinomial model in $(\phi_1,\phi_2)$.
The conditional distribution of $X_1$ given $T=X_1+X_2=t$  is Bi$(t; \phi_1/(1+\phi_1))$, so that, by Table \ref{tab_nef-qvf}, the fiducial density for $\phi_1/(1+\phi_1)$ is Be$(x_1+1/2, t-x_1+1/2)$ and $h^G_{x_1,t}(\f_1)$ coincides with \eqref{eqn:ref_eta}.
 From the marginal distribution of $T$, which  is Bi$(n; \phi_2(1+\phi_1))$, it is possible to derive the fiducial density
 \be
 h^G_{t}(\phi_2|\phi_1)=\frac{\Gamma(n+1)}{\Gamma(t+1/2)\Gamma(n-t+1/2)}
 (1+\phi_1)^{t+1/2}\phi_2^{t-1/2}(1-(1+\phi_1)\phi_2)^{n-t-1/2}
 \ee
  so that the joint fiducial density is  $h^G_{x_1,x_2}(\phi_1,\phi_2)= h_{x_1,t}^G(\phi_1)h^G_{t}(\phi_2|\phi_1)$  which coincides with the joint reference  posterior.

\section{Conclusions and final remarks} \label{sec:conclusions}

We have suggested a way to construct a fiducial distribution  which depends on the  inferential importance ordering of the  parameter components.  Our proposal  appears to be quite simple to apply and, even if it is not so general as the theory suggested by Hannig, has some  advantages in connection with the modern confidence distribution theory and it is strictly related to objective Bayesian analysis.

In  complex models an exact analysis is generally not possible, but approximate results can be derived working with asymptotic  distributions.
 In \citet{Veronese:2014}, starting from the sufficient statistic, an expansion up to the first order of the fiducial distribution for the mean  parameter of a real NEF is provided.
  This  result can be extended to arbitrary regular models starting from the maximum likelihood estimator of the parameter.  When the maximum likelihood estimator is not sufficient a better fiducial distribution can be obtained using an ancillary statistic, as suggested in Section \ref{sec:step-by-step fiducial}. To this aim  the magic formula $p^*$ given by \citet{Barndorff:1983}, which provides an approximation of the conditional distribution of the maximum likelihood estimator given an ancillary statistic, can be fruitfully adopted. Furthermore, these asymptotic results appear to be  strictly connected with the theory of matching priors, i.e.  priors that ensure approximate frequentist validity of posterior credible set. Notice that also these priors crucially depend on the inferential ordering of the parameters, see \citet{Tibshirani:1989} and \citet{Datta-Muk:2004}. However, a normal approximation of the fiducial distribution, when it can be established and  is enough for the analysis, can be proved to be order-invariant. These type of results will be discussed in a forthcoming paper.

\section*{Acknowledgements}

This research was supported by grants from Bocconi University.

\section*{Appendix}
\subsection*{A1: Useful results on cr-NEFs} \label{sec:A1}
Some technical aspects related to cr-NEFs are the following.
\begin{enumerate}
\item
A NEF is a cr-NEF if and only if
the principal $k
\times k$ matrix of the variance function does not depend on
$\mu_{k+1}, \ldots \mu_{d}$, for  $k=1,\ldots,d-1$.
\item The Fisher information matrix relative to the $\gfi$-parametrization is diagonal
with the $kk$-th element depending only on $\gfi_{[k]}$.
%

\item
The cumulant transform  $M_k(\f_k; \gx_{[k-1]})$ of the $k$-th conditional density is given by
\ben  \label{eqn:cond-cum}
M_k(\f_k; \gx_{[k-1]}) & = &
\sum_{j=1}^{k-1}A_{kj}(\f_k)x_j+B_k(\f_k),
\een
for some functions $A_{kj}$ and $B_k$.
\item
The conditional expectation of $X_{k}$ given
$\gX_{[k-1]}=\gx_{[k-1]}$ is linear in $\gx_{[k-1]}$, because it is the
gradient of (\ref{eqn:cond-cum}).
\item
The parameter $\mu_k$ depends on $\gfi$ only through $\gfi_{[k]}$, because from (\ref{eqn:cond-cum})
\ben \label{eqn:mu-phi}
\mu_{k}=  \sum_{j=1}^{k-1}\frac{\partial A_{kj}(\f_k)}{\partial\f_{k}}\mu_j+\frac{\partial B_k(\f_{k})}{\partial \f_{k}}.
\een
\item
Using (\ref{eqn:c_r_r}) and (\ref{eqn:cond-cum}), it can be checked that
\ben \label{eqn:theta_phi}
\te_{k}=\f_{k}-\sum_{u=k+1}^{d}{A}_{uk}(\f_{u}), \quad \mbox{and} \quad M(\gte) = \sum_{k=1}^d B_k(\f_{k}(\gte)).
\een
As a consequence of the first part of (\ref{eqn:theta_phi}), there exists a function  $g_{k}$ such that
 \ben \label{eqn:phi-theta}
\f_{k}= \te_{k} + g_{k}(\te_{k+1}, \ldots, \te_{d}).
 \een
\end{enumerate}
Of course all the previous formulas
hold  for $k=1,\ldots,d$, with the understanding that
components that lose meaning for a specific $k$ are set to zero.
\bigskip

A NEF has a {\em simple quadratic variance
function} (SQVF) if the $ij$-th element of its variance-covariance matrix, seen as a function of the mean parameter $\gmu=(\mu_1,\ldots,\mu_d)$, can be written as
$ V_{ij}(\gmu)=q\mu_i\mu_j+\sum_{k=1}^d \mu_kL_{ij}^{(k)}+C_{ij},
$
where  $q$ is a real constant and $\boldsymbol{L}^{(k)}$, $k=1,\dots,d$ and $\boldsymbol{C}$
are constant $d\times d$ symmetric matrices.
Any NEF-SQVF can be obtained, via a nonsingular affine transformation, from one of the {\em basic} families: Poisson/normal ($q=0$), multinomial ($q=-1/N$, $N$ positive integer), negative-multinomial ($q=1/R$, $R$ positive integer),
ne\-gative-multinomial/gam\-ma/normal ($q=1/R$) and
negative-multinomial/hyperbolic-secant ($q=1/R$), see \citet{Casalis:1996} for a detailed description of these distributions.
The $ij$-th element of the variance function $\gV(\gmu)$ of a basic NEF-SQVF is
\ben \label{basic}
V_{ij}(\gmu)=q\mu_i\mu_j+\sum_{k=1}^d z_{ik} \mu_k+C_{ij},
\een
where $z_{ij}=z_{ji}$,  $i,j \in \{1,\ldots,d\}$, are constants. The values of $z_{ii}$ for the basic NEF-SQVFs, together with other technical details, are given in the proof of Corollary \ref{cor:fid-mu-sqvf}.

\subsection*{A2: Proofs} \label{sec:A3}

\noindent
\emph{Proof of Proposition \ref{prop:invariance}.}\\
By the standard change-of-variable rule  applied to the first integral in \eqref {inv-fid}, it is enough to show that
\ben \label{hfi-hla}
h_\gt^{\gfi}(\gfi(\gla))|J_{\gfi}(\gla)|=h_\gt^{\gla}(\gla),
\een
where $J_{\gfi}(\gla)$ is the Jacobian of the transformation from $\gfi$ to $\gla$.
Now, from \eqref{fid-general} we have
\be
h_{\gt}^{\gfi}(\gfi(\gla))&=&\prod_{k=1}^d h_{\gt_{[k]}, \gt_{-[d]}}(\f_{d-k+1}(\lambda_{[d-k+1]})|\gfi_{[d-k]}(\gla_{[d-k]}))\\
&=& \prod_{k=1}^d \left.\left|\frac{\partial}{\partial \f_{d-k+1}} F_{\f_{d-k+1}}(t_k|\gt_{[k-1]},\gt_{-[d]}; \gfi_{[d-k]})\right|\right|_{\gfi_{[d-k+1]}=\gfi_{[d-k+1]}(\gla_{[d-k+1]})}
\ee
while
\be
J_{\gfi}(\gla)=\prod_{k=1}^d \frac{\partial \f_{d-k+1}(\gla_{[d-k+1]})}{\partial \lambda_{d-k+1}},
\ee
because the transformation from $\gfi$ to $\gla$ is lower triangular.
It follows from the last two formulas and the chain rule that
\be
h_\gt^{\gfi}(\gfi(\gla))|J_{\gfi}(\gla)|=
\prod_{k=1}^d \left|\frac{\partial}{\partial \la_{d-k+1}} F_{\la_{d-k+1}}(t_k|\gt_{[k-1]},\gt_{-[d]}; \gla_{[d-k]})\right|,
\ee
where  $F_{\la_{d-k+1}}$ is the distribution function of $T_k$ given
$(\gT_{[k-1]}=\gt_{[k-1]}, \gT_{[-d]}=\gt_{[-d]})$ in the $\gla$ parameterization. The equality \eqref{hfi-hla} follows by applying  \eqref{fid-general-2} to the  model parameterized by $\gla$.
\st

\medskip

\noindent
\emph{Proof of Proposition \ref{prop: KL-geom-mean}.}\\
 Let $p^G (x)=$ $c^{-1}(p_1(x)p_2(x))^{1/2}$, where $c=\int (p_1(x)p_2(x))^{1/2} d\nu(x)$ is the normalizing constant. Then
\ben \label{kl-geom}
KL(q|p_1)&+&KL(q|p_2)= \int \left(\log \frac{q(x)}{p_1(x)}\right) q(x) d\nu(x)+ \int \left(\log \frac{q(x)}{p_2(x)}\right) q(x) d\nu(x)  \nonumber\\
&=&\int \left(\log \frac{q^2(x)}{p_1(x)p_2(x)}\right) q(x) d\nu(x)
= 2 \int \left(\log \frac{q(x)}{Cp^G(x)}\right) q(x) d\nu(x) \nonumber \\
&=& 2 \int \left(\log \frac{q(x)}{p^G(x)}\right) q(x) d\nu(x)  - 2 \log c
= 2 KL(q|p^G) -2 \log c.
\een
Because $c$ does not depend on $q$, it follows that the functional in (\ref{kl-geom})
achieves its minimum (equal to $-2 \log c$) if and only if $KL(q|p^G)=0$, i.e. $q=p^G$.
\st

\medskip

\noindent
\emph{Proof of Proposition \ref{prop: _geom_mean_between}.}   \\
We only prove that $H_s(\theta) <H_s^G (\theta)$; the other inequality  can be shown in the same way.
Using (\ref{hns}) and  (\ref{left-fid}), we can write
\ben
\frac{h_s^G(\theta)}{h_s(\theta)}&=& \nonumber
\frac{1}{c}\frac{\sqrt{h_s(\theta) h_s^\ell (\theta)}}{h_s (\theta)}
=\frac{1}{c}\frac{\sqrt{h_s(\theta) (h_s(\theta)+\frac{\partial}{\partial\theta}p_\theta(s))}}{h_s(\theta)}\\ \label{decr}
&=&\frac{1}{c}\sqrt{1+\frac{\frac{\partial}{\partial\theta}p_\theta(s)}{h_s(\theta)}}=
\frac{1}{c}\sqrt{1+\gamma_s(\theta)}.
\een
By hypothesis $\gamma_s(\theta)$ is decreasing and thus, from  (\ref{decr}), $h_s^G(\theta)/h_s(\theta)$ is also decreasing on $\Theta$. This is a sufficient condition for $H_s(\theta) <H_s^G (\theta)$, see \citet[Theorem 1.C.1]{Shaked:2007}. \st
\medskip

\noindent
\emph{Proof of Corollary \ref{cor:gamma_s_decreasing_NEF}}.  \\
 Let $p_\theta (s) =\exp(\theta s-n M(\theta)$ be the probability mass function (with respect to a measure $\nu$) of a real NEF, with $\theta$ the natural parameter. Fixing $s \in  {\cal S}_0$, we can write
\ben \label{recipr_gamma}
\gamma_s (\theta) &=&  \left(\frac{\partial p_\theta(s)}{\partial \theta}\right)
/\left(\frac{\partial}{\partial \theta}(1- F_\theta(s))\right)
=\; \frac{(s-nM^\prime (\theta))\exp(\theta s-nM(\theta))}{\sum_{t=s+1}^{+\infty}(t-nM^\prime (\theta)) \exp(\theta t -nM(\theta))} \nonumber \\
&=& \left( \sum_{t=s+1}^{+\infty}\frac{t-nM^\prime (\theta)}{s-nM^\prime (\theta)} \exp\{(t-s)\theta\} \right)^{-1}.
\een
The elements in the sum \eqref{recipr_gamma} are continuous and increasing functions of $\theta$ in both intervals for which $\theta<\widehat{\theta}_s$ and $\theta>\widehat{\theta}_s$, where  $\widehat{\theta}_s=(M^\prime)^{-1}(s/n)$. Thus    $\gamma_s(\theta)$ is  decreasing in these intervals. Moreover, $\gamma_s (\theta)$ is equal to  zero for  $\theta=\widehat{\theta}_s$, positive for $\theta<\widehat{\theta}_s$ and negative for $\theta>\widehat{\theta}_s$, because the denominator of $\gamma_s (\theta)$ is  $h_s(\theta)$, which is positive. Then $\gamma_s(\theta)$ is decreasing over all $\Theta$ and from Proposition \ref{prop: _geom_mean_between} the result follows.  \st

\medskip

\noindent
\emph{Proof of Proposition \ref{prop:comparison_Hg_Ha}.}\\
 In order to prove the proposition, it is sufficient to show that there exist $\theta_1$ and $\theta_2$ in $\Theta$,  $\theta_1<\theta_2$, such that  $h_s^A(\theta_i)=h_s^G(\theta_i)$, $i=1,2$, with $h_s^A(\theta) < h_s^G(\theta)$ for $\theta_1<\theta<\theta_2$ and $h_s^A(\theta) > h_s^G(\theta)$ otherwise, see \citet[proof of Theorem 3.A.44]{Shaked:2007}. Thus we analyze the sign of $h_s^A (\theta)- h_s^G (\theta)$. We can write
 \be
h_s^A (\theta)- h_s^G (\theta) &=& \frac{1}{2}(h_s(\theta)+h_s^\ell(\theta))-\frac{1}{c}\sqrt{h_s(\theta)h_s^\ell(\theta)} \\
  &=& h_s(\theta)\left(\frac{1}{2}(2+\gamma_s(\theta)) -\frac{1}{c} \sqrt{1+\gamma_s(\theta)}\right),
  \ee
  so that the sign of the difference $h_s^A (\theta)- h_s^G (\theta)$ is a function of $\gamma_s(\theta)$ only. First notice that, by a standard property of the arithmetic and geometric means, $c=\int \sqrt{h_s(\theta)h_s^\ell(\theta)}d\theta < \int (h_s(\theta)+h_s^\ell(\theta))/2 \hspace{1 mm} d\theta=1$ for all $\theta$.
  After some straightforward algebra, it can be seen that  $h_s^A (\theta)- h_s^G (\theta)=0$ when (and only when) $\gamma_s (\theta)=2c^{-2}((1-c^2)-\sqrt{1-c^2})=k_1$ or $\gamma_s (\theta)=2c^{-2}((1-c^2)+\sqrt{1-c^2})=k_2$, with  $k_1 \in (-1,0)$ and $k_2>0$.  Moreover we have    $h_s^A (\theta)<h_s^G (\theta)$ for $k_1<\gamma_s(\theta)<k_2$ and $h_s^A (\theta)>h_s^G (\theta)$ for $\gamma_s(\theta)<k_1$ or $\gamma_s(\theta)>k_2$.  By assumption, $\gamma_s(\theta)$ is decreasing on $\Theta$ from $+\infty$ to -1, so that there exist $\theta_1$ and $\theta_2$, with $\gamma_s(\theta_1)=k_2$ and $\gamma_s(\theta_2)=k_1$, satisfying the sufficient condition stated at the beginning of the proof.  \st


\medskip

\noindent
\emph{Proof of Proposition \ref{prop:sufficiency}.}\\
First notice that if we use for constructing the fiducial distribution the sufficient statistic $S$, which has the same dimension of the parameter,  we do not need an ancillary statistic.
Furthermore, $\gT$ is a one-to-one transformation of  $\gX$, and thus  $\gS$ is a function of $\gT=(\gTqd,\gTmqd)$ but, since  $\gS$ is  complete, it is stochastically independent  of $\gTmqd$ by Basu's Theorem. As a consequence,   $\gS_{[k]}$ is also independent of  $\gTmqd$ and thus
\ben \label{suff}
& &\mbox{Pr}_{\f_{d-k+1}}(S_k\leq s_k \mid\gS_{[k-1]}=
\gs_{[k-1]}; \gfi_{[d-k]}) = \nonumber\\
& &\mbox{Pr}_{\f_{d-k+1}}(S_k\leq s_k \mid\gS_{[k-1]}=\gs_{[k-1]},\gT_{-[d]}=\gt_{-[d]}; \gfi_{[d-k]}).
\een
From the one-to-one lower triangular transformation $\gs=\textbf{g}(\gt_{[d]}, \gt_{-[d]})$, we have that $s_k=g_k(t_k,\gt_{[k-1]},\gt_{-[d]})$, with $g_k$ invertible with respect to $t_k$, so that (assuming $g_k$ increasing) \eqref{suff} becomes
\be
& &\mbox{Pr}_{\f_{d-k+1}}(g_k(T_k,\gT_{[k-1]},\gTmqd)\leq s_k \mid\gT_{[k-1]}=\gt_{[k-1]},\gT_{-[d]}=\gt_{-[d]}; \gfi_{[d-k]})=\\
& &\mbox{Pr}_{\f_{d-k+1}}(T_k \leq g_k^{-1}(s_k,\gT_{[k-1]},\gTmqd) \mid\gT_{[k-1]}=\gt_{[k-1]},\gT_{-[d]}=\gt_{-[d]}; \gfi_{[d-k]})=\\
& &\mbox{Pr}_{\f_{d-k+1}}(T_k \leq t_k \mid\gT_{[k-1]}=\gt_{[k-1]},\gT_{-[d]}=\gt_{-[d]}; \gfi_{[d-k]}),
\ee
which proves the proposition.  \st

\medskip

\noindent
\emph{Proof of Proposition \ref{prop:fiducial-cr-NEF}.}\\
Because each conditional distribution of $X_{k}$ given $\gX_{[k-1]}=\gx_{[k-1]}$ belongs to a  NEF with natural parameter $\f_k$, using \eqref{Hns} we have that $\Hgsk(\f_k)$ is a distribution function for  $\f_k$. The result follows from the postulated independence among the $\f_k$'s.
\st

\medskip

\noindent
\emph{Proof of Proposition \ref{prop:fid-NRF-SQVF}.}\\
Formulas \eqref{basic-fid} and \eqref{fiducial-NM/G/N} derive by a direct application of \eqref{fid-NEF-SQVF} to the conditional distributions of the different families. For a detailed description  of the cr-NEFs  involved,  see \citet[proof of Theorem 3]{Consonni:2001}.
\st

\medskip

\noindent
\emph{Proof of Corollary \ref{cor:fid-mu-sqvf}.}\\
First notice that the fiducial distribution $h_{\gs}(\gmu)$ can be more easily obtained via a double transformation, namely
\be
h_{\gs}^{\gmu}(\gmu)=h_{\gs}^{\gte}(\gte(\gmu))|J_{\gfi}(\gte(\gmu))||J_{\gte}(\gmu)|,
\ee
where the Jacobian   $|J_{\gfi}(\gte)|=1$  for \eqref{eqn:phi-theta}, and
\ben \label{J-teta-phi}
|J_{\gte}(\gmu)|\propto
\mbox{det}\{{\gV}(\gmu)\}^{-1}= \exp\left\{-\sum_{k=1}^d \te_k(\gmu)z_k -q(d+1) M(\gte(\gmu)) \right\},
\een
 see \citet[pag. 34]{Gutierrez:Smith:1997} for the proportionality  relationship and \citet[Prop. 1]{Consonni:2004} for the  equality. We consider now each family.

 \medskip

\noindent
\emph{Poisson/normal family} with $m$ Poisson components. We have $q=0$,  $\f_k=\te_k$;   $z_k=1$,  $\te_k=\log(\mu_k)$ and $B_k(\f_k)=\exp(\f_k)$ for $k=1,\dots,m$,  while  $z_k=0$, $\te_k= \mu_k/\s^2$ and $B_k(\f_k)=\s^2 \f_k^2/2$ for $k=m+1,\dots,d$,    where $\s^2$ is the known variance of the normal components.
Then from \eqref{J-teta-phi}, it  follows that $|J_{\gfi}(\gte(\gmu))|=\prod_{k=1}^m \mu_k^{-1}$,  and thus using
 \eqref{basic-fid}, the result \eqref{fid-poisson-mu} follows.

\medskip

\noindent
\emph{Multinomial family}. Using the relationships  in Example 1, \eqref{J-teta-phi} gives
$|J_{\gfi}(\gte(\gmu))|=(N-\sum_{k=1}^d \mu_k)^{-1}\prod_{k=1}^d \mu_k^{-1}$ and using \eqref{fid-general} we obtain \eqref{fid-multinomial mu}.

\medskip

\noindent
\emph{Negative-multinomial family}. We have $q=1/R$, $R>0$,  $z_k=1$, $\f_k(\gte)=\te_k-\log(1-\sum_{u=k+1}^d e^{\te_u})$,  $\te_k=\log(\mu_k) -\log(R+\sum_{j=1}^d \mu_j)$, and  $B_k(\f_k)=-R\log(1-e^{\f_k})$ for $k=1,\dots,d$.
Then from \eqref{J-teta-phi}, it  follows that $|J_{\gfi}(\gte(\gmu))|=(R+\sum_{k=1}^d \mu_k)^{-1}\prod_{k=1}^d \mu_k^{-1}$,  and thus using
 \eqref{basic-fid}, the result \eqref{fid-neg-multinomial-mu} follows.

\medskip

\noindent
\emph{Negative-multinomial/gamma/normal family (with an $m$ dimensional negative-multinomial component)}. We have $q=1/R$, $R>0$;  $z_k=1$ and $\f_k=\log(\mu_k/(R+\sum_{j=1}^k \mu_j))$, $k=1,\dots,m$;  $z_{m+1}=0$ and $\f_{m+1}=-(R+\sum_{j=1}^m \mu_j)/\mu_{m+1}$; $z_k=0$ and $\f_k=\mu_k/\mu_{m+1}$, $k=m+2,\dots,d$. In this case it is convenient to compute the fiducial density of $\gmu$ directly from  \eqref{fiducial-NM/G/N}. Observing that the Jacobian of the transformation from $\gfi$ to $\gmu$ is
\be \label{J-phi-mu-NM-gamma-normal}
|J_{\gfi}(\gmu)| =  \left(R+\sum_{k=1}^m \mu_k\right)^{-1} \left(\prod_{k=1}^m \mu_k^{-1}\right)  \left(R+\sum_{j=1}^m \mu_j\right)\mu^{-2}_{m+1} \, \mu_{m+1}^{-d+m+1}
\ee
 and using the previous expression of $\f_k$, the density \eqref{fid-neg-multinomial-gamma-normal-mu} follows.
 \st

\noindent
\emph{Proof of Proposition \ref{prop:fid location or scale}.}\\
Let $\gX$ be an  i.i.d. sample of size $n$, with $X_i \sim p_{\te}(x_i)=f(x_i-\te)$, i.e. $\te$ is a location parameter, and consider the transformation $Z_1=X_1$, $Z_i=X_i-X_1$, $i=2, \ldots, n$, whose Jacobian is one. Then, setting $\gz=(z_1, \dots, z_n)$
\be
H_{\gx}(\te)=H_{\gz}(\te)=1-F_{\te}(z_1|z_2, \dots, z_n)=\frac{\int_{z_1}^{+\infty} f(t-\te) \prod_{i=2}^n f(t+z_i-\te) dt}{\int_{-\infty}^{+\infty} f(t-\te) \prod_{i=2}^n f(t+z_i-\te) dt}.
\ee
Using now the substitution  $m=-t+\te+z_1$ in the previous two integrals, and recalling  that $\pi^J(\te) \propto 1$ we obtain
\be
\frac{\int^{\te}_{-\infty} f(z_1-m) \prod_{i=2}^n f(z_1+z_i-m) dm}{\int_{-\infty}^{+\infty} f(z_1-m)
\prod_{i=2}^n f(z_1+z_i-m) dm} &=& \frac{\int^{\te}_{-\infty}  \prod_{i=1}^n f(x_i-m) \pi^J(m) dm}{\int_{-\infty}^{+\infty} \prod_{i=1}^n f(x_1-w) \pi^J(m) dm} \\
& =&  \int^{\te}_{-\infty}\pi^J(m|\gx) dm.
 \ee
The result relative to the scale parameter follows recalling that the model $p_{\te}(x)=f(x/\te)/\te$
can be transformed in a model with location parameter $\mu$ setting $y=\log(x)$ and $\mu=\log(\te)$.
In this case a constant prior on $\mu$ is equivalent to a prior on $\te$ proportional to $1/\te$.
\st

\medskip

\noindent

\noindent
\emph{Proof of Proposition \ref{prop:fid location-scale}.}\\
Let $\gX$ be an  i.i.d. sample of size $n$, with $X_i \sim p_{\te,\s}(x_i)=f((x_i-\te)/\s)/\s$, $i=1,\ldots,n$ and notice that the absolute value of the Jacobian of the transformation from  $\gx$ to $(x_1,z_2,\gz)$, with $z_2=x_2-x_1$, $\gz=(z_3,\dots, z_n)$ and $z_j=(x_j-x_1)/z_2$, $j=3,\dots,n$, is $|z_2|^{n-2}$. Furthermore, the reference prior $\pi^R(\te,\s)$ is order-invariant and can be written as $\pi^R(\te|\s) \pi^R(\s)$, where  $\pi^R(\te|\s) \propto 1$ and $\pi(\s) \propto 1/\s$, see \citet{Steel:1999}. Working conditionally on $\s$ we can thus apply Proposition \ref{prop:fid location or scale}  to conclude that the reference posterior and the fiducial distribution for $\te$ given $\s$ coincide. It remains to show that
\ben \label{post-sigma}
\pi^R(\s|\gx) & = &\pi^R(\s|x_1,z_2,\gz) \propto  \int_{-\infty}^{\infty} p_{\,\te,\s}(x_1,z_2,\gz) \pi^R(\te|\s) \pi^R(\s) d\te \\
&\propto & \int_{-\infty}^{\infty}\frac{1}{\s^{n+1}}
f\left(\frac{x_1-\te}{\s}\right) f\left(\frac{z_2+x_1-\te}{\s}\right) \prod_{i=3}^nf\left(\frac{z_2z_i+x_1-\te}{\s}\right) d\te \nonumber
\een
corresponds to the fiducial density $h_{z_2,\gz}(\s)$.

We have
\ben \label{fid-sigma}
H_{z_2,\gz}(\s)& = & 1-F_{\s}(z_2|\gz)=\int_{z_2}^{+\infty}\int_{-\infty}^{+\infty} p_{\,\te,\s}^{X_1,Z_2|\gZ}(t,w|\gz)/p^{\gZ}(\gz) dt dw \\
& = &\frac{\int_{z_2}^{+\infty} \int_{-\infty}^{+\infty}
\frac{|w|^{n-2}}{\s^{n}}
f\left(\frac{t-\te}{\s}\right) f\left(\frac{w+t-\te}{\s}\right) \prod_{i=3}^nf\left(\frac{wz_i+t-\te}{\s}\right) dt dw}{p^{\gZ}(\gz)}, \nonumber
\een
where the density of $p^{\gZ}(\gz)$ does not depend on the parameters because $\gz$ is ancillary.
Assuming $z_2>0$, and using the transformation $m= x_1-v(t-\te)/\s, v=z_2\s/w$,  which implies  $t=\s(x_1-m)/v + \te, w=z_2\s/v$, with Jacobian $z_2\s^2/v^3$, the fiducial distribution $H_{z_2,\gz}(\s)$ in \eqref{fid-sigma} becomes
\be
\frac{\int_{0}^{\s} \int_{-\infty}^{+\infty}
\frac{z_2^{n-1}}{v^{n+1}}
f\left(\frac{x_1-m}{v}\right) f\left(\frac{z_2+x_1-m}{v}\right) \prod_{i=3}^nf\left(\frac{z_2z_i+x_1-m}{v}\right) dm dv}{p^{\gZ}(\gz)}
.
\ee
Taking the derivative with respect to $\s$, it is immediate to see that the fiducial density for $\s$ coincides with the posterior distribution given in \eqref{post-sigma}.

If $z_2<0$, and applying to the integral the same transformation  used in the previous case, we have
\be
   F_{\s}(z_2|\gz) &=&  \frac{\int^{z_2}_{-\infty} \int_{-\infty}^{+\infty}
\frac{|w|^{n-2}}{\s^{n}}
f\left(\frac{t-\te}{\s}\right) f\left(\frac{w+t-\te}{\s}\right) \prod_{i=3}^nf\left(\frac{wz_i+t-\te}{\s}\right) dt dw}{{p^{\gZ}(\gz)}}
 \\
&=& -\frac{\int_{0}^{\s} \int_{-\infty}^{+\infty}
\frac{(-z_2)^{n-1}}{v^{n+1}}
f\left(\frac{x_1-m}{v}\right) f\left(\frac{z_2+x_1-m}{v}\right) \prod_{i=3}^nf\left(\frac{z_2z_i+x_1-m}{v}\right) dm dv}{p^{\gZ}(\gz)}
,
\ee
so that again the derivative with respect to $\s$ of $H_{z_2,\gz}(\s) = 1-F_{\s}(z_2|\gz)$ leads to \eqref{post-sigma}.
\st

\medskip

The following lemma will be used in the proof of Proposition \ref{prop: fid-prior-cr-NEF}.
\bl \label{lemma: jeffreys-reference}
Consider a cr-NEF on $\Rex^d$, with the $k$-th diagonal element in the Fisher information matrix given by $I_{kk}(\gfi)=a_k(\f_k)b_k(\gfi_{[k-1]})$. Then the $d$-group (order-invariant) reference prior $\pi^R$ for $\gfi=(\f_1, \ldots, \f_d)$ is
\ben \label{ref_jef}
 \pi^R(\gfi)=\prod_{k=1}^d\pi^J_k(\f_k)\propto \prod_{k=1}^d (a_k(\f_k))^{1/2},
\een
where $\pi^J_k(\f_k)$ is the Jeffreys prior  obtained from the conditional distribution of $X_k$ given $\gX_{[k-1]}=\gx_{[k-1]}$.
\el

\noindent
\emph{Proof of Lemma \ref{lemma: jeffreys-reference}}\\
First observe that $\gmu_{[k]}$ is a one-to-one transformation of $\gfi_{[k]}$ and that the information matrix  $\gI(\gfi)$ of a cr-NEF is diagonal,  see  Appendix A1 (points 2 and  5). From  \eqref{eqn:c_r_r}  and \eqref{eqn:cond-cum}, the $kk$-th element of $\gI$ is
\be
 I_{kk}(\gfi_{[k]})& =& -E^{\gX}_{\gfi}\left(\frac{\partial^2}{\partial \f_k^2}\log p_{\f_{k}}(x_k|\gx_{[k-1]}; \f_k)\right)= -E^{\gX}_{\gfi}(M''_k(\f_k; \gx_{[k-1]}))\\
 &=& \sum_{j=1}^{k-1}A''_{kj}(\f_k)\mu_j(\gfi_{[j]})+B''_k(\f_k).
\ee
 Under the assumption in the proposition,   we can write
$I_{kk}(\gfi_{[k]})=a_k(\f_k)b^*_k(\gmu_{[k-1]}(\gfi_{[k-1]}))$. From \citet{Datta:1995}, it follows that  the reference prior   on $\gfi$ is order-invariant and is given by the last product in \eqref{ref_jef}.

Consider now the Jeffreys prior on $\f_k$ obtained from  $p_{\f_{k}}(x_k|\gx_{[k-1]})$. This is proportional to the square root of
\be
-E^{X_k|\gx_{[k-1]}}(M''_k(\f_k; \gx_{[k-1]}))=
 \sum_{j=1}^{k-1}A''_{kj}(\f_k)x_j+B''_k(\f_k)= a_k(\f_k)b^*_k(\gx_{[k-1]}),
 \ee
where again the last equality holds by the assumption in the proposition. Thus the product of the $d$ Jeffreys priors is equal to \eqref{ref_jef} and the result holds.
\st

\noindent
\emph{Proof of Proposition \ref{prop: fid-prior-cr-NEF}.}\\
Due to the independence of the $\f_k$'s,  a fiducial prior for  $\gfi$ exists if and only if there exists a
fiducial prior for each $\f_k$.  Because the conditional distribution of $S_k$ given $\gS_{[k-1]}=\gs_{[k-1]}$
belongs to a real NEF with natural parameter $\f_k$, the result of the first part of the proposition follows from Proposition \ref{prop:FID_JEF}.

 The first statement of the second part of the proposition follows checking directly  the form of the conditional distributions of the basic NEF-SQVFs and using again Proposition \ref{prop:FID_JEF}. The second statement follows from the remark stated before the proposition and from Lemma \ref{lemma: jeffreys-reference}.
\st


\begin{thebibliography}{100}
%
\bibitem[{Barndorff-Nielsen(1983)}]{Barndorff:1983}
{Barndorff-Nielsen, O.} (1983).
\newblock {On a formula for the distribution of the maximum likelihood estimator}.
\newblock \textit{Biometrika} \textbf{70}, 343--365.

%
\bibitem[{Berger(2006)}]{Berger:2006}
{Berger, J.~O.}   (2006).
\newblock The case for objective Bayesian analysis.
\newblock \textit{Bayesian Analysis} \textbf{1}, 385--402.

\bibitem[{Berger \& Bernardo (1992)}]{Berger:1992}
{Berger, J.~O.} \& {Bernardo,  J.~M.}
  (1992).
\newblock Ordered group reference priors with application to a multinomial
problem.
\newblock \textit{Biometrika} \textbf{79}, 25--37.


\bibitem[{Berger et~al.(2015)}]{Berger:2015}
{Berger, J.~O.} {Bernardo,  J.~M.} \& {Sun, D.}
  (2015).
\newblock Overall objective priors.
\newblock \textit{Bayesian Analysis} \textbf{10}, 189--221.

\bibitem[{Bernardo(1979)}]{Bernardo:1979}
{Bernardo, J.~M.}   (1979).
\newblock Reference posterior distributions for Bayesian inference.
\newblock \textit{J. R. Stat. Soc. Ser. {\rm B}} \textbf{41}, 113--147.

\bibitem[{Bernardo \& Ramon(1998)}]{Bernardo:1998}
{Bernardo, J.M.}  \& {Ramon, J.M.}
  (1998).
\newblock  An introduction to Bayesian reference analysis: inference on the ratio of multinomial parameters.
\newblock \textit{The Statistician} \textbf{28}, 101--135.

\bibitem[{Bernardo \& Smith(1994)}]{Bernardo:1994}
{Bernardo, J.~M.} \& {Smith, A.~F.~M.} (1994).
\newblock \textit{{ Bayesian Theory}}.
\newblock  Wiley: Chichester.
%

\bibitem[{Casalis(1996)}]{Casalis:1996}
{Casalis, M.} (1996).
\newblock {The 2d+4 simple quadratic natural exponential families on $\Rex^d$}.
\newblock \textit{Ann. Statist.} \textbf{24}, 1828--1854.

\bibitem[{Consonni \& Veronese(2001)}]{Consonni:2001}
{Consonni, G.}  \& {Veronese, P.}
  (2001).
\newblock Conditionally reducible natural exponential families and enriched conjugate priors.
\newblock \textit{Scand. J. Stat.} \textbf{28}, 377--406.

\bibitem[{Consonni et~al.(2004)}]{Consonni:2004}
{Consonni, G.}, {Veronese, P.} \& {Guti\'errez-Pe\~na, E.}
  (2004).
\newblock Reference priors for exponential families with simple quadratic variance function.
\newblock \textit{J. Multivariate Anal.} \textbf{88}, 335--364.


\bibitem[{Crowder \& Sweeting(1989)}]{Crowder:1989}
{Crowder, M.}, \& {Sweeting, T.} (1989).
\newblock {Bayesian inference for a bivariate binomial distribution}.
\newblock \textit{Biometrika} \textbf{76}, 599--603.

%
\bibitem[{Datta \& Ghosh(1995)}]{Datta:1995}
{Datta, G.~S.} \& {Ghosh, M.} (1995).
\newblock {Some Remarks on Noninformative Priors.}
\newblock \textit{J.  Amer.
Statist. Assoc.} \textbf{90},  1357--1363.
%
\bibitem[{Datta \& Ghosh(1996)}]{Datta:1996}
{Datta, G.~S.} \& {Ghosh, M.} (1996).
\newblock {On the Invariance of Noninformative Priors.}
\newblock \textit{Ann. Statist.} \textbf{24}, 141--159.




\bibitem[{Datta \& Mukerjee(2004)}]{Datta-Muk:2004}
{Datta, G.~S.} \& {Mukerjee, R.} (2004).
\newblock \textit{{ Probability matching priors: higher order asymptotics (Lecture Notes in Statistics)}}.
\newblock  Springer: New York.

\bibitem[{Dawid \& Stone(1982)}]{Dawid:1982}
{Dawid, A.~P.}  \& {Stone, M.} (1982).
\newblock {The functional-model basis of fiducial inference.}
\newblock \textit{Ann. Statist.} \textbf{10}, 1054--1074.

\bibitem[{Dempster(1963)}]{Dempster:1963}
{Dempster, A.~P.}  (1963).
\newblock {Further examples of inconsistencies in the fiducial argument.}
\newblock \textit{Ann. Statist.} \textbf{34}, 884--891.

\bibitem[{Fern\'andez \& Steel (1999)}]{Steel:1999}
{Fern\'andez, C.} \& {Steel, F.~J.} (1999).
\newblock {Reference priors for the general location-scale model}.
\newblock \textit{Statist. Prob. Lett.} \textbf{43}, 377--384.




\bibitem[{Fisher(1930)}]{Fisher:1930}
{Fisher, R.~A.} (1930).
\newblock {Inverse probability.}
\newblock \textit{Proceedings of the Cambridge Philosophical Society}
 \textbf{26}, 4, 528--535.

\bibitem[{Fisher(1935)}]{Fisher:1935}
{Fisher, R.~A.} (1935).
\newblock {The fiducial argument in statistical inference}.
\newblock \textit{Ann. Eugenics} \textbf{VI}, 91--98.

\bibitem[{Fisher(1973)}]{Fisher:1973}
{Fisher, R.~A.} (1973).
\newblock \textit{{Statistical methods and scientific inference}}.
\newblock Hafner Press: New York.

\bibitem[{Fraser(1961)}]{Fraser:1961}
{Fraser, D.~A.~S.} (1961).
\newblock {On fiducial  inference}.
\newblock \textit{Ann. Math. Statist.} \textbf{32}, 661--676.




%

\bibitem[{Guti\'errez-Pe\~na \& Smith(1997)}]{Gutierrez:Smith:1997}
{Guti\'errez-Pe\~na, E.}  \& {Smith, A.~F.~M.}
  (1997).
\newblock Exponential and Bayesian conjugate families: review and extensions (with discussion).
\newblock \textit{Test} \textbf{6}, 1--90.

\bibitem[{Hannig(2009)}]{Hannig:2009}
{Hannig, J.}  (2009).
\newblock On generalized fiducial inference.
\newblock \textit{Statist. Sinica} \textbf{19}, 491--544.

\bibitem[{Hannig(2013)}]{Hannig:2013}
{Hannig, J.}  (2013).
\newblock Generalized fiducial inference via discretization.
\newblock \textit{Statist. Sinica} \textbf{23}, 489--514.

\bibitem[{Hannig \& Iyer(2008)}]{Hannig:2008}
{Hannig, J.}  \& {Iyer, H.}  (2008).
\newblock Fiducial intervals for variance components in an unbalanced two-component normal mixed linear model.
\newblock \textit{J. Amer. Statist. Assoc.} \textbf{103}, 854--865.

\bibitem[{Hannig et~al.(2007)}]{Hannig:2007}
{Hannig, J.}, {Iyer, H.~K.} \& {Wang, C.~M.}
  (2007).
\newblock Fiducial approach to uncertainty assessment accounting for error due to instrument resolution.
\newblock \textit{Metrologia} \textbf{44}, 476--483.

\bibitem[{Hannig et~al.(2016)}]{Hannig:2016}
{Hannig, J.}, {Iyer, H.~K.},  {Lai, R.~C.~S.} \& {Lee T.~C.~M.}
 (2016).
\newblock Generalized Fiducial Inference: A Review and New Results.
\newblock \textit{J. American Statist. Assoc.} \textbf{44}, 476--483.


%
%


\bibitem[{Johnson et ~al.(2005)}]{Johnson:2005}
{Johnson, N.~L., Kemp, W.~A. \& Kotz, J.~P.} (2005).
\newblock \textit{{Univariate discrete distributions}}.
\newblock Wiley: New York.



\bibitem[{Krishnamoorthy \& Lee(2010)}]{Krishnamoorthy:2010}
{Krishnamoorthy, K.} \& {Lee, M.}
  (2010).
\newblock  Inference for functions of parameters in discrete distributions based on fiducial approach: binomial and Poisson cases.
\newblock \textit{J.  Statist. Plann.  Inference.} \textbf{140}, 1182--1192.


\bibitem[{Lehmann \& Romano(2005)}]{Lehmann:2005}
{Lehmann, E.~L. \& Romano, J.~P.} (2005).
\newblock \textit{{Testing statistical hypotheses}}.
\newblock Springer: New York.


\bibitem[{Lindley(1958)}]{Lindley:1958}
{Lindley, D.~V.} (1958).
\newblock {Fiducial distributions and Bayes theorem}.
\newblock \textit{J. R. Stat. Soc. Ser. {\rm B}} \textbf{20}, 102--107.
%
%
\bibitem[{Martin \& Liu(2013)}]{Martin:2013}
{Martin, R.}  \& {Liu, C.}
  (2013).
\newblock Inferential models: a framework for prior-free posterior probabilistic inference.
\newblock \textit{J. Amer. Statist. Assoc.} \textbf{108}, 301--313.

%
%
\bibitem[{Petrone \& Veronese(2010)}]{Petrone:2010}
{Petrone, S.}  \& {Veronese, P.}
  (2010).
\newblock Feller operators and mixture priors in Bayesian nonparametrics.
\newblock \textit{Statist. Sinica} \textbf{20}, 379--404.





\bibitem[{Polson \& Wasserman(1990)}]{Polson:1990}
{Polson, N.}, \& {Wasserman, L.} (1990).
\newblock {Prior distributions for the bivariate binomial}.
\newblock \textit{Biometrika} \textbf{77}, 901--904.



\bibitem[{Schweder \& Hjort(2002)}]{Schweder:2002}
{Schweder, T.} \& {Hjort, N.~L.}
  (2002).
\newblock Confidence and likelihood.
\newblock \textit{Scand. J. Stat.} \textbf{29}, 309--332.



\bibitem[{Schweder \& Hjort(2016)}]{Schweder:2016}
{Schweder, T.} \& {Hjort, N.~L.} (2016).
\newblock  \textit{Confidence, likelihood and probability}.
\newblock London: Cambridge University Press.




\bibitem[{Shaked \& Shanthikumar(2007)}]{Shaked:2007}
{Shaked, M. \& Shanthikumar, J.~G.} (2007).
\newblock \textit{{Stochastic orders}}.
\newblock Springer: New York.

\bibitem[{Singh et~al.(2005)}]{Singh:2005}
{Singh, K.}, {Xie, M.}  \& {Strawderman, M.}
  (2005).
\newblock Combining information through confidence distribution.
\newblock \textit{Ann. Statist.} \textbf{33}, 159--183.

\bibitem[{Stein(1959)}]{Stein:1959}
{Stein, C.} (1959).
\newblock \textit{An example of wide discrepancy between fiducial and confidence intervals}.
\newblock \textit{Ann. Math. Statist.} \textbf{30}, 877--880.


%
\bibitem[{Taraldsen \& Lindqvist(2013)}]{Taraldsen:2013}
{Taraldsen, G.}  \& {Lindqvist, B.~H.}
  (2013).
\newblock Fiducial theory and optimal inference.
\newblock \textit{Ann. Statist.} \textbf{41}, 323--341.

\bibitem[{Tibshirani(1989)}]{Tibshirani:1989}
{Tibshirani, R.}
  (1989).
\newblock Noninformative priors for one parameter of many.
\newblock \textit{Biometrika} \textbf{76}, 604--608.



\bibitem[{V\&M(2015)}]{Veronese:2014}
{Veronese, P.}  \& {Melilli, E.}
  (2015).
\newblock Fiducial and confidence distributions for real exponential families.
\newblock \textit{Scand. J. Stat.} \textbf{42}, 471--484.


%






%
%
\bibitem[{Wandler \& Hannig(2012)}]{Wandler:2012}
{Wandler, D.}  \& {Hannig, J.}
  (2012).
\newblock A fiducial approach to multiple comparisons.
\newblock \textit{J. Statist. Plann. Inference} \textbf{142}, 878--895.


%


\bibitem[{Wilkinson(1977)}]{Wilkinson:1977}
{Wilkinson, G.~N.} (1977).
\newblock {On resolving the controversy in statistical inference}.
\newblock \textit{J. R. Stat. Soc. Ser. {\rm B}} \textbf{39}, 119--171.




\end{thebibliography}
\end{document}